\documentclass[a4paper,12pt]{article}
\usepackage{amssymb}
\usepackage{latexsym}
\usepackage{color}
\newcommand{\comment}[1]{}  

\newtheorem{thm}{Theorem:}[section]

\newtheorem{cor}[thm]{Corollary:}

\newtheorem{lem}[thm]{Lemma:}
\newtheorem{deff}[thm]{Definition:}

\newtheorem{rem}[thm]{Remark:}

\newtheorem{examples}[thm]{Examples:}

\newcommand{\rk}{\mbox{rank}}

\newcommand{\CZF}{{\mathbf{CZF}}}
\newcommand{\ZFC}{{\mathbf{ZFC}}}
\newcommand{\ZF}{{\mathbf{ZF}}}

\newcommand{\AC}{{\mathbf{AC}}}
\newcommand{\PA}{{\mathbf{PA}}}

\newcommand{\prf}{{\bf Proof: }}

\newcommand{\gdw}{\,\leftrightarrow\,}

\newcommand{\beqs}{\begin{eqnarray*}}
\newcommand{\eeqs}{\end{eqnarray*}}

\newcommand{\BI}{{\mathbf{BI}}}

\newcommand{\beq}{\begin{eqnarray}}
\newcommand{\eeq}{\end{eqnarray}}

\newcommand{\IZF}{{\mathbf{IZF}}}

\newcommand{\suc}{{\mathrm{suc}}}

\newcommand{\N}{{\mathbb N}}

\newcommand{\PAs}[1]
{#1^+}

\newcommand{\bes}{\begin{eqnarray*}}
\newcommand{\ees}{\end{eqnarray*}}

\newcommand{\DP}{{\mathbf{DP}}}

\newcommand{\bluff}{ \blue\bf }

 \newcommand{\TS}{\,\Rightarrow\,}
  \newcommand{\TSI}{\,\Rightarrow\,} 

\newcommand{\LL}{{\mathcal{L}}}

\newcommand{\TI}{{\mathrm {TI}}}

\newcommand{\GCH}{{\mathbf{GCH}}}

\newcommand{\bgf}{\begin{frame}}
\newcommand{\efr}{\end{frame}}

\newcommand{\ee}{\end{itemize}\end{frame}}

\newcommand{\blue}{\color{blue}}
\newcommand{\red}{\color{red}}

\newcommand{\provxii}[4]{\provx{I_{_{\!\infty}}\!\!#1}{#2}{#3}{#4}}

\newcommand{\vieh}{\varphi}

\newcommand{\reddish}{\em }

\def\fCenter{\ \Rightarrow\ }

\def\BPmessage{Proof Tree (bussproofs) style macros. Version 0.6c.}
\def\EnableBpAbbreviations{%
	\let\AX\Axiom
	\let\AXC\AxiomC
	\let\UI\UnaryInf
	\let\UIC\UnaryInfC
	\let\BI\BinaryInf
	\let\BIC\BinaryInfC
	\let\TI\TrinaryInf
	\let\TIC\TrinaryInfC
	\let\LL\LeftLabel
	\let\RL\RightLabel
	\let\DP\DisplayProof
}

\def\ScoreOverhang{4pt}			
\def\ScoreOverhangLeft{\ScoreOverhang}
\def\ScoreOverhangRight{\ScoreOverhang}

\def\extraVskip{2pt}			
\def\ruleScoreFiller{\hrule}		

\def\defaultScoreFiller{\ruleScoreFiller}  
\def\defaultBuildScore{\buildSingleScore}  

\def\defaultHypSeparation{\hskip.2in}   

\def\labelSpacing{3pt}		

\def\proofSkipAmount{\vskip.8ex plus.8ex minus.4ex}

\def\theHypSeparation{\defaultHypSeparation}
\def\alwaysScoreFiller{\defaultScoreFiller}	
\def\alwaysBuildScore{\defaultBuildScore}
\def\theScoreFiller{\alwaysScoreFiller}	
\def\buildScore{\alwaysBuildScore}   
\def\hypKernAmt{0pt}	

\def\defaultLeftLabel{}
\def\defaultRightLabel{}

\def\myTrue{Y}
\def\bottomAlignFlag{N}
\def\centerAlignFlag{N}

\expandafter\ifx\csname newenvironment\endcsname\relax%
\def\makeatletter{\catcode`\@=11\relax}
\def\makeatother{\catcode`\@=12\relax}
\makeatletter
\def\newcount{\alloc@0\count\countdef\insc@unt}
\def\newdimen{\alloc@1\dimen\dimendef\insc@unt}
\def\newskip{\alloc@2\skip\skipdef\insc@unt}
\def\newbox{\alloc@4\box\chardef\insc@unt}
\makeatother
\else
\newenvironment{prooftree}%
{\begin{center}\proofSkipAmount \leavevmode}%
{\DisplayProof \proofSkipAmount \end{center} }
\fi

\def\thecur#1{\csname#1\number\theLevel\endcsname}

\newcount\theLevel    
\global\theLevel=0		
\newcount\myMaxLevel
\global\myMaxLevel=0
\newbox\myBoxA		
\newbox\myBoxB
\newbox\myBoxC
\newbox\myBoxD
\newbox\myBoxLL		
\newbox\myBoxRL
\newdimen\thisAboveSkip		
\newdimen\thisBelowSkip		
\newdimen\newScoreStart		
\newdimen\newScoreEnd
\newdimen\newCenter
\newdimen\displace
\newdimen\leftLowerAmt
\newdimen\rightLowerAmt
\newdimen\scoreHeight
\newdimen\scoreDepth

\setbox\myBoxLL=\hbox{\defaultLeftLabel}%
\setbox\myBoxRL=\hbox{\defaultRightLabel}%

\def\allocatemore{%
	\ifnum\theLevel>\myMaxLevel%
		\expandafter\newbox\curBox%
		\expandafter\newdimen\curScoreStart%
		\expandafter\newdimen\curCenter%
		\expandafter\newdimen\curScoreEnd%
		\global\advance\myMaxLevel by1%
	\fi%
}

\def\prepAxiom{%
	\advance\theLevel by1%
	\edef\curBox{\thecur{myBox}}%
	\edef\curScoreStart{\thecur{myScoreStart}}%
	\edef\curCenter{\thecur{myCenter}}%
	\edef\curScoreEnd{\thecur{myScoreEnd}}%
	\allocatemore%
}

\def\Axiom$#1\fCenter#2${%
	\prepAxiom%
	\setbox\myBoxA=\hbox{$\mathord{#1}\fCenter\mathord{\relax}$}%
	\setbox\myBoxB=\hbox{$#2$}%
	\global\setbox\curBox=%
	     \hbox{\hskip\ScoreOverhangLeft\relax%
		\unhcopy\myBoxA\unhcopy\myBoxB\hskip\ScoreOverhangRight\relax}%
	\global\curScoreStart=0pt \relax
	\global\curScoreEnd=\wd\curBox \relax
	\global\curCenter=\wd\myBoxA \relax
	\global\advance \curCenter by \ScoreOverhangLeft%
	\ignorespaces
}

\def\AxiomC#1{		
	\prepAxiom%
	\setbox\myBoxA=\hbox{#1}%
	\global\setbox\curBox =%
		\hbox{\hskip\ScoreOverhangLeft\relax%
                        \unhcopy\myBoxA\hskip\ScoreOverhangRight\relax}%
        \global\curScoreStart=0pt \relax
        \global\curScoreEnd=\wd\curBox \relax
        \global\curCenter=.5\wd\curBox \relax
        \global\advance \curCenter by \ScoreOverhangLeft%
	\ignorespaces
}

\def\prepUnary{%
	\ifnum \theLevel<1
		\errmessage{Hypotheses missing!}
	\fi%
	\edef\curBox{\thecur{myBox}}%
	\edef\curScoreStart{\thecur{myScoreStart}}%
	\edef\curCenter{\thecur{myCenter}}%
	\edef\curScoreEnd{\thecur{myScoreEnd}}%
}

\def\UnaryInf$#1\fCenter#2${%
	\prepUnary%
	\buildConclusion{#1}{#2}%
	\joinUnary%
	\resetInferenceDefaults%
	\ignorespaces%
}

\def\UnaryInfC#1{
	\prepUnary%
	\buildConclusionC{#1}%
	\joinUnary%
	\resetInferenceDefaults%
	\ignorespaces%
}

\def\prepBinary{%
	\ifnum\theLevel<2
		\errmessage{Hypotheses missing!}
	\fi%
	\edef\rcurBox{\thecur{myBox}}
	\edef\rcurScoreStart{\thecur{myScoreStart}}%
	\edef\rcurCenter{\thecur{myCenter}}%
	\edef\rcurScoreEnd{\thecur{myScoreEnd}}%
	\advance\theLevel by-1
	\edef\lcurBox{\thecur{myBox}}
	\edef\lcurScoreStart{\thecur{myScoreStart}}%
	\edef\lcurCenter{\thecur{myCenter}}%
	\edef\lcurScoreEnd{\thecur{myScoreEnd}}%
}

\def\BinaryInf$#1\fCenter#2${%
	\prepBinary%
	\buildConclusion{#1}{#2}%
	\joinBinary%
	\resetInferenceDefaults%
	\ignorespaces%
}

\def\BinaryInfC#1{%
	\prepBinary%
	\buildConclusionC{#1}%
	\joinBinary%
	\resetInferenceDefaults%
	\ignorespaces%
}

\def\prepTrinary{%
	\ifnum\theLevel<3
		\errmessage{Hypotheses missing!}
	\fi%
	\edef\rcurBox{\thecur{myBox}}
	\edef\rcurScoreStart{\thecur{myScoreStart}}%
	\edef\rcurCenter{\thecur{myCenter}}%
	\edef\rcurScoreEnd{\thecur{myScoreEnd}}%
	\advance\theLevel by-1
	\edef\ccurBox{\thecur{myBox}}
	\edef\ccurScoreStart{\thecur{myScoreStart}}%
	\edef\ccurCenter{\thecur{myCenter}}%
	\edef\ccurScoreEnd{\thecur{myScoreEnd}}%
	\advance\theLevel by-1
	\edef\lcurBox{\thecur{myBox}}
	\edef\lcurScoreStart{\thecur{myScoreStart}}%
	\edef\lcurCenter{\thecur{myCenter}}%
	\edef\lcurScoreEnd{\thecur{myScoreEnd}}%
}

\def\TrinaryInf$#1\fCenter#2${%
	\prepTrinary%
	\buildConclusion{#1}{#2}%
	\joinTrinary%
	\resetInferenceDefaults%
	\ignorespaces%
}

\def\TrinaryInfC#1{%
	\prepTrinary%
	\buildConclusionC{#1}%
	\joinTrinary%
	\resetInferenceDefaults%
	\ignorespaces%
}

\def\buildConclusion#1#2{
        \setbox\myBoxA=\hbox{$\mathord{#1}\fCenter\mathord{\relax}$}%
        \setbox\myBoxB=\hbox{$#2$}%
	\setbox\myBoxC =%
	      \hbox{\hskip\ScoreOverhangLeft\relax%
		\unhcopy\myBoxA\unhcopy\myBoxB\hskip\ScoreOverhangRight\relax}%
	\newScoreStart=0pt \relax%
	\newCenter=\wd\myBoxA \relax%
	\advance \newCenter by \ScoreOverhangLeft%
	\newScoreEnd=\wd\myBoxC%
}

\def\buildConclusionC#1{
	\setbox\myBoxA=\hbox{#1}%
	\setbox\myBoxC =%
		\hbox{\hbox{\hskip\ScoreOverhangLeft\relax%
                        \unhcopy\myBoxA\hskip\ScoreOverhangRight\relax}}%
	\newScoreStart=0pt \relax%
	\newCenter=.5\wd\myBoxC \relax%
	\newScoreEnd=\wd\myBoxC%
        \advance \newCenter by \ScoreOverhangLeft%
}

\def\joinUnary{
	\global\advance\curCenter by -\hypKernAmt%
	\ifnum\curCenter<\newCenter%
		\displace=\newCenter%
		\advance \displace by -\curCenter%
		\kernUpperBox%
	\else%
		\displace=\curCenter%
		\advance \displace by -\newCenter%
		\kernLowerBox%
	\fi%
        \ifnum \newScoreStart < \curScoreStart %
		\global \curScoreStart = \newScoreStart \fi%
	\ifnum \curScoreEnd < \newScoreEnd %
		\global \curScoreEnd = \newScoreEnd \fi%
	\ifnum \curScoreStart<\wd\myBoxLL%
		\global\displace = \wd\myBoxLL%
		\global\advance\displace by -\curScoreStart%
		\kernUpperBox%
		\kernLowerBox%
	\fi%
	\buildScore%
	\buildScoreLabels%
	\global \setbox \curBox =%
		\vbox{\box\curBox%
			\vskip\thisAboveSkip \relax%
			\nointerlineskip\box\myBoxD%
			\vskip\thisBelowSkip \relax%
			\nointerlineskip\box\myBoxC}%
	\global \curScoreStart=\newScoreStart%
	\global \curScoreEnd=\newScoreEnd%
	\global \curCenter=\newCenter%
}
\def\kernUpperBox{%
		\global\setbox\curBox =%
			\hbox{\hskip\displace\box\curBox}%
		\global\advance \curScoreStart by \displace%
		\global\advance \curScoreEnd by \displace%
		\global\advance\curCenter by \displace%
}
\def\kernLowerBox{%
		\global\setbox\myBoxC =%
			\hbox{\hskip\displace\unhbox\myBoxC}%
		\global\advance \newScoreStart by \displace%
		\global\advance \newScoreEnd by \displace%
		\global\advance\newCenter by \displace%
}

\def\joinBinary{
	\setbox\myBoxA=\hbox{\theHypSeparation}
	\lcurScoreEnd=\rcurScoreEnd%
	\advance\lcurScoreEnd by\wd\lcurBox%
	\advance\lcurScoreEnd by\wd\myBoxA%
	\displace=\lcurScoreEnd%
	\advance\displace by -\lcurScoreStart%
	\lcurCenter=.5\displace%
	\advance\lcurCenter by\lcurScoreStart%
	\setbox\lcurBox=%
		\hbox{\box\lcurBox\unhcopy\myBoxA\box\rcurBox}%
	\displace=\newCenter%
	\advance\displace by -.5\newScoreStart%
	\advance\displace by -.5\newScoreEnd%
	\advance\lcurCenter by \displace%
	\edef\curBox{\lcurBox}%
	\edef\curScoreStart{\lcurScoreStart}%
	\edef\curScoreEnd{\lcurScoreEnd}%
	\edef\curCenter{\lcurCenter}%
	\joinUnary%
}

\def\joinTrinary{
	\setbox\myBoxA=\hbox{\theHypSeparation}
	\lcurScoreEnd=\rcurScoreEnd%
	\advance\lcurScoreEnd by\wd\lcurBox%
	\advance\lcurScoreEnd by\wd\ccurBox%
	\advance\lcurScoreEnd by2\wd\myBoxA%
	\displace=\lcurScoreEnd%
	\advance\displace by -\lcurScoreStart%
	\lcurCenter=.5\displace%
	\advance\lcurCenter by\lcurScoreStart%
	\setbox\lcurBox=%
		\hbox{\box\lcurBox\unhcopy\myBoxA\box\ccurBox%
				  \unhcopy\myBoxA\box\rcurBox}%
	\displace=\newCenter%
	\advance\displace by -.5\newScoreStart%
	\advance\displace by -.5\newScoreEnd%
	\advance\lcurCenter by \displace%
	\edef\curBox{\lcurBox}%
	\edef\curScoreStart{\lcurScoreStart}%
	\edef\curScoreEnd{\lcurScoreEnd}%
	\edef\curCenter{\lcurCenter}%
	\joinUnary%
}

\def\DisplayProof{%
	\ifnum \theLevel=1 \relax \else
		\errmessage{Proof tree badly specified.}%
	\fi%
	\edef\curBox{\thecur{myBox}}%
	\ifx\bottomAlignFlag\myTrue%
		\displace=0pt%
	\else%
		\displace=.5\ht\curBox%
		\ifx\centerAlignFlag\myTrue\relax
		\else%
		      	\advance\displace by -3pt%
		\fi%
	\fi%
	\leavevmode%
	\lower\displace\hbox{\copy\curBox}%
	\global\theLevel=0%
	\global\def\alwaysBuildScore{\defaultBuildScore}
	\global\def\alwaysScoreFiller{\defaultScoreFiller}
	\def\bottomAlignFlag{N}
	\def\centerAlignFlag{N}
	\resetInferenceDefaults%
	\ignorespaces
}

\def\buildSingleScore{
	\displace=\curScoreEnd%
	\advance \displace by -\curScoreStart%
	\global\setbox \myBoxD =%
		\hbox to \displace{\expandafter\xleaders\theScoreFiller\hfill}%
}

\def\buildDoubleScore{
	\buildSingleScore%
	\global\setbox\myBoxD=%
		\hbox{\hbox to0pt{\copy\myBoxD\hss}\raise2pt\copy\myBoxD}%
}

\def\buildNoScore{
	\global\setbox\myBoxD=\hbox{\vbox{\vskip1pt}}%
}

\def\doubleLine{%
	\gdef\buildScore{\buildDoubleScore}
	\ignorespaces
}

\def\LeftLabel#1{%
	\global\setbox\myBoxLL=\hbox{{#1}\hskip\labelSpacing}%
	\ignorespaces
}
\def\RightLabel#1{%
	\global\setbox\myBoxRL=\hbox{\hskip\labelSpacing #1}%
	\ignorespaces
}

\def\buildScoreLabels{%
	\scoreHeight = \ht\myBoxD%
	\scoreDepth = \dp\myBoxD%
	\leftLowerAmt=\ht\myBoxLL%
	\advance \leftLowerAmt by -\dp\myBoxLL%
	\advance \leftLowerAmt by -\scoreHeight%
	\advance \leftLowerAmt by \scoreDepth%
	\leftLowerAmt=.5\leftLowerAmt%
	\rightLowerAmt=\ht\myBoxRL%
	\advance \rightLowerAmt by -\dp\myBoxRL%
	\advance \rightLowerAmt by -\scoreHeight%
	\advance \rightLowerAmt by \scoreDepth%
	\rightLowerAmt=.5\rightLowerAmt%
	\displace = \curScoreStart%
	\advance\displace by -\wd\myBoxLL%
	\global\setbox\myBoxD =%
		\hbox{\hskip\displace%
			\lower\leftLowerAmt\copy\myBoxLL%
			\box\myBoxD%
			\lower\rightLowerAmt\copy\myBoxRL}%
	\global\thisAboveSkip = \ht\myBoxLL%
	\global\advance \thisAboveSkip by -\leftLowerAmt%
	\global\advance \thisAboveSkip by -\scoreHeight%
	\ifnum \thisAboveSkip<0 %
		\global\thisAboveSkip=0pt%
	\fi%
	\displace = \ht\myBoxRL%
	\advance \displace by -\rightLowerAmt%
	\advance \displace by -\scoreHeight%
	\ifnum \displace<0 %
		\displace=0pt%
	\fi%
	\ifnum \displace>\thisAboveSkip %
		\global\thisAboveSkip=\displace%
	\fi%
	\global\thisBelowSkip = \dp\myBoxLL%
	\global\advance\thisBelowSkip by \leftLowerAmt%
	\global\advance\thisBelowSkip by -\scoreDepth%
	\ifnum\thisBelowSkip<0 %
		\global\thisBelowSkip = 0pt%
	\fi%
	\displace = \dp\myBoxLL%
	\advance\displace by \rightLowerAmt%
	\advance\displace by -\scoreDepth%
	\ifnum\displace<0 %
		\displace = 0pt%
	\fi%
	\ifnum\displace>\thisBelowSkip%
		\global\thisBelowSkip = \displace%
	\fi
	\global\thisAboveSkip = -\thisAboveSkip%
	\global\thisBelowSkip = -\thisBelowSkip%
	\global\advance\thisAboveSkip by\extraVskip
	\global\advance\thisBelowSkip by\extraVskip
}

\def\resetInferenceDefaults{%
	\global\def\theHypSeparation{\defaultHypSeparation}%
	\global\setbox\myBoxLL=\hbox{\defaultLeftLabel}%
	\global\setbox\myBoxRL=\hbox{\defaultRightLabel}%
	\global\def\buildScore{\alwaysBuildScore}%
	\global\def\theScoreFiller{\alwaysScoreFiller}%
	\gdef\hypKernAmt{0pt}
}

\def\provx#1#2#3#4{
\setbox1=\hbox{\kern1.5pt$\scriptstyle#3$}
\def\zeichen{#2}
\ifx\zeichen\empty\setbox0=\hbox to .75em{}\else\setbox0=\hbox
{\kern1.5pt$\scriptstyle#2$}\fi
\dimen1=\dp0 \ifdim \dimen1=0pt
\advance \dimen1 by 1.5ex \else \advance \dimen1 by 1.2ex
\fi\dimen3=2ex\dimen4=.5ex\ifdim \wd0<\wd1 \dimen2=\wd1 \else \dimen2=\wd0
\fi\hbox{$#1\hskip 5pt minus5pt\vrule height\dimen3
depth\dimen4\raise\dimen1\copy0\hskip-1\wd0 \lower\ht1
\copy1\hskip-1\wd1\vrule width\dimen2 height.7ex depth-.6ex\hskip3pt
minus1.5pt#4\hskip2pt plus2pt minus2pt$}}

\def\prov#1#2#3{
\setbox1=\hbox{\kern1.5pt$\scriptstyle#2$}
\def\zeichen{#1}
\ifx\zeichen\empty\setbox0=\hbox to .75em{}\else\setbox0=\hbox
{\kern1.5pt$\scriptstyle#1$}\fi
\dimen1=\dp0
\ifdim \dimen1=0pt
\advance \dimen1 by 1.5ex \else \advance \dimen1 by 1.2ex
\fi\dimen3=2ex\dimen4=.5ex\ifdim \wd0<\wd1 \dimen2=\wd1 \else \dimen2=\wd0
\fi\hbox{\hskip0pt plus 4pt
$\vrule height\dimen3
depth\dimen4\raise\dimen1\copy0\hskip-1\wd0
\lower\ht1\copy1\hskip-1\wd1\vrule width\dimen2 height.7ex depth-.6ex
\hskip3pt minus1.5pt#3\hskip2pt plus2pt minus2pt$}}

\def\prv#1#2{
\setbox1=\hbox{\kern1.5pt$\scriptstyle#2$}
\ifx\zeichen\empty\setbox0=\hbox to .75em{}\else\setbox0=\hbox
{\kern1.5pt$\scriptstyle#1$}\fi
\dimen1=\dp0 \ifdim \dimen1=0pt
\advance \dimen1 by 1.5ex \else \advance \dimen1 by 1.2ex
\fi\dimen3=2ex\dimen4=.5ex\ifdim \wd0<\wd1 \dimen2=\wd1 \else \dimen2=\wd0
\fi\hbox{\hskip.5em$\vrule height\dimen3
depth\dimen4\raise\dimen1\copy0\hskip-1\wd0
\lower\ht1\copy1\hskip-1\wd1\vrule width\dimen2 height.7ex depth-.6ex
\hskip3pt minus1.5pt$}}

\mathchardef\str='1066
\def\negprov#1#2#3{
\setbox1=\hbox{\kern1.5pt$\scriptstyle#2$}
\setbox4=\hbox{$\str$}
\def\zeichen{#1}
\ifx\zeichen\empty\setbox0=\hbox to 1em{}\else\setbox0=\hbox
{\kern1.5pt$\scriptstyle#1$}\fi
\dimen1=\dp0
\ifdim \dimen1=0pt
\advance \dimen1 by 1.5ex \else \advance \dimen1 by 1.2ex
\fi\dimen3=2ex\dimen4=.5ex\ifdim \wd0<\wd1 \dimen2=\wd1 \else \dimen2=\wd0
\fi\hbox{\hskip.5em$\kern-1.9pt\raise1pt\copy4\kern-\wd4\kern1.9pt\vrule height\dimen3
depth\dimen4\raise\dimen1\copy0\hskip-1\wd0
\lower\ht1\copy1\hskip-1\wd1\vrule width\dimen2 height.7ex depth-.6ex
\hskip3pt minus1.5pt#3\hskip2pt plus2pt minus2pt$}}

\def\goed#1{\setbox5=\hbox{$#1$}\dimen1=.25em \dimen2=\dimen1 \advance \dimen2
by -1pt\hbox{\raise.65\ht5 \hbox{\vrule height.5\ht5 depth0pt width.4pt\vrule
height.5\ht5 width\dimen1 depth-.48\ht5}\kern-\dimen2\copy5\kern-\dimen2
\raise.65\ht5 \hbox{\vrule height .5\ht5 width\dimen1 depth-.48\ht5\vrule
height.5\ht5 depth 0pt width.4pt}\hskip4pt plus2pt minus2pt}}

\def\mod#1#2{
\def\zeichen{#1}
\hbox{\hskip 2pt plus3pt minus 2pt\vrule width.5pt height2ex depth.5ex
\vbox{\ifx\zeichen\empty\hbox to .75em{}\else
\hbox{\kern1.5pt $\scriptstyle#1$}\fi
\kern2pt
\hrule
\kern1.7pt
\hrule\kern1.7pt}
\hskip3pt minus 2pt$#2$}\hskip2pt
plus3pt minus2pt}

\def\notmod#1#2{\hbox{\hskip 2pt plus 3pt minus 3pt\vrule width.5pt
height2ex depth.5ex
\vbox{\hbox{\kern1.5pt $\scriptstyle#1$}\kern3pt
\setbox0=\hbox{\kern2pt$\scriptstyle/$}
\hrule
\kern-1.7pt
\copy0
\kern-\ht0
\kern 1.7pt
\hrule\kern1.7pt}n
\hskip3pt minus 2pt$#2$}\hskip2pt
plus3pt minus2pt}
\def\sq{\hbox{\rlap{$\sqcap$}$\sqcup$}}
\def\qed{\ifmmode\sq\else{\unskip\nobreak\hfil\penalty50\hskip1em\null
\nobreak\hfil\sq\parfillskip=0pt\finalhyphendemerits=0\endgraf}\fi\medskip}

\def\lleq{\hbox{\hskip3pt minus3pt\kern1pt\lower4pt
\vbox{\hbox{$\scriptstyle\ll$}
\kern-7pt\hbox{\kern1pt$\scriptstyle=$}}\hskip3pt minus 3pt}}

\mathchardef\res='1152
\mathchardef\qin='1062
\mathchardef\qprec='1036
\mathchardef\qless='474
\mathchardef\dpkt='72

\newcommand{\On}{{\mathrm{ON}}}
\newcommand{\DDD}{{\mathcal D}}
\newcommand{\ITT}{{\mathcal{IL}}}
\begin{document}

\title{Remarks on Barr's theorem: Proofs in geometric theories}
\author{Michael Rathjen\\
{\small Department of Pure Mathematics, University of Leeds}\\
{\small Leeds LS2 9JT, United Kingdom}\\
 {\tt \scriptsize
rathjen@maths.leeds.ac.uk}}

\maketitle
\begin{abstract} A theorem, usually  attributed to Barr, yields that (A) geometric implications deduced in classical ${\mathcal L}_{\infty\omega}$ logic from geometric theories also have intuitionistic proofs.  Barr's theorem is of a  topos-theoretic nature and its proof is non-constructive. In the literature
one also finds mysterious comments about the capacity of this theorem to remove the axiom of choice from derivations.
  This  article  investigates the proof-theoretic side of Barr's theorem and also aims to shed some light on the axiom of choice  part. More concretely, a constructive proof of
   the Hauptsatz for ${\mathcal L}_{\infty\omega}$ is given and is put to use to arrive at a simple proof of (A) that is formalizable in constructive set theory and Martin-L\"of type theory.

\end{abstract}

\section{Introduction}
  A {\em signature} $\Sigma$ consists of constant symbols, function symbols, and relation symbols together with an assignment of a unique positive integer (arity) to
any object of the latter two kinds.
A language
${\mathcal L}$ is comprised of a signature $\Sigma$ and formation rules, i.e., rules for forming formulae over $\Sigma$.
The familiar Tarskian way of assigning meaning to the symbols of $\mathcal L$ proceeds by associating set-theoretic objects to them, notably functions and relations construed set-theoretically, giving rise
to the notion of (set-theoretic) {\em structure for $\mathcal L$} and {\em model
of $T$} for any theory $T$ in the language of $\mathcal L$.
There is also a more general notion of structure in a sufficiently rich category $\mathcal C$. For example, if $\mathcal C$ has finite products, then any equational language (i.e., equality being the sole relation symbol and equations the only formulae)
allows for interpretation in $\mathcal C$, by viewing terms as morphisms
and function symbol application as composition.
Another prominent example is the interpretation of the typed  $\lambda$-calculus in cartesian closed categories (cf. \cite{lambek-scott}). For still richer languages one must impose more conditions on $\mathcal C$. If one wants to extend this idea to higher order logic, then $\mathcal C$
is required to be a topos (cf. \cite[D1.2]{johnstone2}).
This extra level of generality of interpretation, however, comes with a penalty to pay in that
 only intuitionistically valid consequences can be guaranteed to survive the interpretation.

 One is often interested in transferring results from the category of sets,
 $\mathbf{Set}$, where classical logic, the axiom of choice and more reign,
 to an arbitrary topos $\mathcal E$.  This is possible, for instance, for the following (non first-order) assertion:
 \begin{center}{\em
 All modules over fields are flat.\footnote{This is just a simple example. Flatness of a module $M$, a notion introduced by Serre in 1956, is usually defined
 by saying that tensoring with $M$ preserves injectivity. An equivalent way of expressing in ${\mathcal L}_{\omega_1\omega}$ that $M$ is a flat $R$-module for a ring $R$  is the following: For all $m\in\mathbb N$,
 whenever $x_1,\ldots,x_m\in M$ and $r_1,\ldots,r_m\in R$ satisfy $\sum r_ix_i=0$, then there exist
 $y_1,\ldots,y_n\in M$ and $a_{ij}\in R$ such that $x_i=\sum a_{ij}y_j$ and $\sum r_ia_{ij}=0$.}}
 \end{center}

 A result that ensures this transfer is commonly called Barr's Theorem (see e.g.\cite[p.515]{mm})\footnote{Disclaimer applying to the entire paper: This is not a paper on the history of certain pieces of mathematics. The attribution of results to persons is borrowed from standard text books or articles in the area, and therefore may
 well be historically inaccurate, as is so often the case.}
 but it can also be
 inferred from cut elimination for the infinitary logic ${\mathcal L}_{\omega_1\omega}$ (see later parts of this paper).
 For this to work, however, the formalization of mathematical notions is important. They have to be chosen carefully, as familiar equivalences are liable to fail in an intuitionistic setting.
 Moreover, to ensure survival of statements it will be important to develop mathematics within (classical) {\em geometric theories} and to couch statements as geometric implications.
 The topos-theoretic result alluded to above is the following.

 \begin{thm}\label{I1} For every Grothendieck topos $\mathcal E$ there exists a complete Boolean algebra $\mathbf B$ and a surjective geometric morphism $\mathrm{Sh}({\mathbf B})\to {\mathcal E}$.
 Here $\mathrm{Sh}({\mathbf B})$ is the topos of sheaves on the Boolean algebra with the usual sup topology. Moreover,
 $\mathrm{Sh}({\mathbf B})$  is a Boolean topos and satisfies the axiom of choice, in the sense that for any epi $e:Y\twoheadrightarrow X$ there exists $s:Y\to X$ such that $e\circ s=1_Y$.
 \end{thm}
 As a consequence\footnote{There are several steps and further theorems involved; cf. \cite[3.1.16]{johnstone2}.}  one arrives at the following insight.
 \begin{cor} If $T$ is a geometric theory and $A$ is a geometric statement
 deducible from $T$ with classical logic, then $A$ is also deducible from
 $T$ with intuitionistic logic, where by logic we mean infinitary
 ${\mathcal L}_{\infty\omega}$-logic.
 \end{cor}
 Though  this Corollary also follows from a syntactic cut elimination
 result for ${\mathcal L}_{\infty\omega}$ (see section \ref{Hauptsatz}), Barr's theorem is often alleged
 to achieve more in that it also allows to eliminate uses of the axiom of choice. This is borne out by the following quotes:
 \begin{quote}
  ``METATHEOREM. If a geometric sentence is deducible from a geometric theory in
classical logic, with the axiom of choice, then it is also deducible from it intuitionistically." G.C. Wraith: {\em Intuitionistic Algebra: Some Recent Developments
in Topos Theory.} Proceedings of the International Congress of Mathematicians,
Helsinki, 1978 331--337.
\end{quote}

\begin{quote}  ``This has the advantage that all such toposes satisfy the Axiom of Choice;
 so we obtain a further conservativity result ..., asserting that uses of the Axiom of Choice may be eliminated from any {\em derivation} of a geometric sequent from geometric hypothesis." P. Johnstone: {\em Sketches of an elephant}, vol. 2, p. 899.
 \end{quote}
 Judging from conversations with logicians and discussions on internet forums (e.g. MathOverflow),
it is probably fair to say that the main appeal of Barr's theorem stems
from its mysterious power to utilize $\AC$ and then subsequently get rid of it.
But can it really perform these wonders?
As a backcloth for the discussion it might be useful to recall some famous
$\AC$-removal results.

\begin{thm}\label{Gloria} Below $\GCH$ stands for the generalized continuum hypothesis. \begin{itemize}
\item[(i)] (G\"odel 1938--1940) If $A$ is a number-theoretic statement and
$\ZFC+\GCH\vdash A$ then $\ZF\vdash A$.
\item[(ii)] (Shoenfield 1961, Platek, Kripke, Silver 1969) If $B$ is a $\Pi^1_4$-statement of second order arithmetic and $\ZFC+\GCH\vdash B$ then $\ZF\vdash B$.
    \item[(iii)] (Goodman 1976, 1978) If $A$ is a number-theoretic statement
    and ${\mathbf{HA}}^{\omega}+\AC_{type}\vdash A$, then ${\mathbf{HA}}\vdash A$. Here ${\mathbf{HA}}$ stands for intuitionistic arithmetic also known as Heyting arithmetic. ${\mathbf{HA}}^{\omega}$ denotes Heyting arithmetic in all finite types with $\AC_{type}$ standing for the collection of all higher type versions
     $\AC_{\sigma\tau}$ of the axiom of choice with $\sigma$, $\tau$ arbitrary finite types.
    \end{itemize}
    \end{thm}
    So should Barr's theorem be added to this list of renowned theorems with $\AC$-eliminatory powers?
The above quotes by Wraith and Johnstone seem to suggest that the addition
of the axiom of choice to a geometric theory does not produce new geometric theorems. But one immediately faces the question of what it means to add $\AC$ to a theory $T$. Here it might be useful to introduce a rough distinction which
separates two ways of doing this. 
The first route, which consists in expressing $\AC$ in the same language as $T$,  will be referred to as an {\em internal} addition of $\AC$. If, on the other hand, $\AC$ is expressed in a richer language with a new sort of objects where the choice functions live; we shall term it an {\em external} addition.
If $T$ is a first-order theory, then adding $\AC$ internally to $T$ requires the language of $T$ to be sufficiently rich. Moreover, internal  $\AC$ forces the choice functions to be objects falling under the first order quantifiers of $T$, and thus, in general, the axioms of $T$ will ``interact"
with $\AC$ in this augmentation.
By contrast, an external addition of $\AC$ refers to a potentially larger universe, where the choice functions needn't be denizens of the realm that the original theory $T$ speaks about. The foregoing distinction is still very coarse, though. For instance the choice functions might be external but they can certainly act on the original objects of $T$. Therefore if one also demands principles of $T$ (e.g. induction) to hold for terms that involve these choice functions (like in the Goodman result) it is conceivable that conservativity will be lost (as is the case with the classical version $\PA^{\omega}+\AC_{type}$ of
${\mathbf{HA}}^{\omega}+\AC_{type}$). Notwithstanding that there are multifarious possibilities to add $\AC$, 
labeling some of them as internal and others as external augmentations provides a useful, if crude, heuristics.
%

In view of Theorem \ref{Gloria} one can also ask if Barr's Theorem can be
beefed up to include more than just $\AC$. For instance, how about the continuum hypothesis, $V=L$, $\lozenge$ and other axioms?

\section{Geometric and $\infty$-geometric theories}
Below we will work in the extension ${\mathcal L}_{\infty\omega}$ of first order logic (${\mathcal L}_{\omega\omega}$) which has all the formulae engendering rules of the latter but also allows  to form infinitely long conjunctions $\bigwedge \Phi$ and disjunctions $\bigvee \Phi$ from any set
$\Phi$ of already constructed formulae. A particulary well-behaved fragment of
${\mathcal L}_{\infty\omega}$ is ${\mathcal L}_{\omega_1\omega}$ where the set $\Phi$ in $\bigwedge\Phi$ and $\bigvee \Phi$ is always required to be countable. Infinitary logics
began to play an important role in logic in the 1950s.\footnote{
``Yet infinitary logic has a long prehistory. Infinitely long formulas were introduced
by C.S. Peirce in the 1880s, used by Schr\"oder in the 1890s, developed
further by L\"owenheim and Lewis in the 1910s, explored by Ramsey and Skolem
in the 1920s, extended by Zermelo and Helmer in the 1930s, studied by Carnap,
Novikov, and Bochvar (and explicitly rejected by G\"odel) in the 1940s,
and exploited by A. Robinson (1951)." \cite{moore}}
\subsection{Geometric theories}
\begin{deff}\label{geometric}{\em The {\em geometric formulae} are inductively defined as follows:
 Every atom is a geometric formula. If $A$, $B$, and $C(a)$ are geometric formulae then so are
 $A\vee B$, $A\wedge B$ and $\exists x\,C(x)$ (where $x$ does not occur in $C(a)$.

Another way of saying this is that a formula is geometric iff it does not contain any of the particles $\to,\neg,\forall$.

 A formula is called a {\em geometric implication} if it is of either form
 $\forall \vec x\,A$ or $\forall \vec x\,\neg A$ or $\forall \vec x\,(A\to B)$
 with $A$ and $B$ being geometric formulae. Here $\forall \vec x$ may be empty.
 In particular geometric formulae and their negations are geometric implications.

 A theory is {\em geometric} if all its axioms are geometric implications.
 }\end{deff}

Below we shall give several examples of geometric theories.\footnote{For more detailed descriptions of these theories and also the ones considered in Section \ref{infge}, see e.g. \cite{bah}, \cite[1.4]{chang-keisler}, \cite[Appendix: Examples]{hodges}.}
 \begin{examples}\label{3.11}{\em
 \begin{itemize}
 \item[(i)] 1. Robinson arithmetic. The language has a constant $0$, a unary successor
function $\suc$ and binary functions $+$ and $\cdot$.
Axioms are the equality axioms and the universal closures of the following.
\begin{enumerate}
\item $\neg \suc(a) = 0$.
\item $\suc(a) = \suc(b)\to a = b$.
\item $a = 0\,\vee\,\exists y \,a = \suc(y)$.
\item $ a + 0 = a$.
\item $ a + \suc(b) = \suc(a + b)$.
\item $a \cdot 0 = 0$.
\item $ a \cdot \suc(b) = a \cdot b + a$
\end{enumerate}
A classically equivalent axiomatization is obtained if (3) is replaced by
$$ \neg a = 0 \to \exists y\, a = \suc(y)$$ but this is not a geometric implication.

\item[(ii)] The theories of  groups, rings, local rings and division rings
have geometric axiomatizations. Local rings are commutative rings with $0\ne 1$ having just one maximal ideal. On the face of it,
the latter property appears to be second order but it can be rendered geometrically as follows:
$$\forall x\,(\exists y\,x\cdot y=1\;\vee\;\exists y\,(1-x)\cdot y=1).$$

\item[(iii)] The theories of fields, ordered fields, algebraically closed fields and
real closed fields have geometric axiomatizations.
To express invertibility of non-zero elements one uses $\forall x\,(x=0\,\vee\,\exists y\,x\cdot y=1)$ rather
than the non-geometric axiom $\forall x\,(x\ne 0\to \exists y\,x\cdot y=1)$.

To express algebraic closure replace axioms

$$s\ne 0\,\to \,\exists x\,sx^n+t_1x^{n-1}+\ldots +t_{n-1}x+t_n=0$$
by
$$s= 0\,\vee \,\exists x\,sx^n+t_1x^{n-1}+\ldots +t_{n-1}x+t_n=0$$
where $sx^k$ is short for $s\cdot x\cdot\ldots \cdot x$ with $k$ many $x$.
\\[1ex]
 Also the theory of {\em differential fields} has a geometric axiomatization.
This theory is written in the language of rings with an additional unary function symbol $\delta$. The axioms
are the field axioms plus $\forall x\forall y\,\delta (x+y)=\delta (x)+\delta (y)$ and
$\forall x\forall y\,\delta (x\cdot y)=x\cdot\delta (y)+y\cdot\delta (x)$.

\item[(iv)] The theory of projective geometry has a geometric axiomatization.

\item[(v)] The theories of equivalence relations, dense linear orders, infinite sets
and graphs also have geometric axiomatizations.

\end{itemize}

 }\end{examples}

\subsection{The infinite geometric case}\label{infge}
Infinitary logics are much more expressive and it is interesting to investigate notions of geometricity in these
 expanded settings.  The infinitary languages we have in mind are such that they accommodate infinite disjunctions $\bigvee \Phi$ and conjunctions $\bigwedge\Phi$, where $\Phi$ is set of (infinitary) formulae.\footnote{It will be assumed that the total number of free variables occurring in the formulae of $\Phi$ is finite. The reason for this commonly found restriction appears to be that in this language  only finitely many variables can be quantified at a time. So if one allowed infinitely many free variables there would be formulae which cannot be closed.}
 This language is customarily denoted by ${\mathcal L}_{\infty\omega}$.

 In this richer language a formula is said to be {\em infinite geometric},
 notated $\infty$-geometric,  if in addition to
 $\vee,\wedge,\exists$ one also allows
 infinite disjunctions $\bigvee \Phi$, where $\Phi$ is already a set of $\infty$-geometric formulae satisfying the above proviso on the number of variables.

 An example  of an axiom expressible in this richer language via a $\infty$-geometric implication
 is the Archimedean axiom:
 $$\forall x \,(x<1\,\vee\,x<1+1\,\vee\ldots\vee\,x<1+\ldots +1\,\vee\ldots)$$
 or in more compact way with $\N^+=\{n\in\N\mid n>0\}$: $$\forall x\,\bigvee_{n\in\N^+} x<n.$$

One often only considers the sublanguage ${\mathcal L}_{\omega_1\omega}$ of ${\mathcal L}_{\infty\omega}$ where the formation of $\bigvee \Phi$ and $\bigwedge \Phi$ is only permissible for countable sets of formulae $\Phi$.

\begin{deff}\label{def-coherent}{\em The {\em $\infty$-geometric formulae} are inductively defined as follows:
 Every atom is a $\infty$-geometric formula. If $A$ and $B$ are $\infty$-geometric formulae then so are
 $A\vee B$ and $A\wedge B$. If $C(a)$ is a $\infty$-geometric formula with all occurrences of $a$ indicated and
 $x$ is a bound variable that does not occur in $C(a)$ then  $\exists x\,C(x)$
 is a $\infty$-geometric formula.
If $\Phi$ is a set of $\infty$-geometric formulae having a finite number of free variables  then
 $\bigvee\Phi$ is a $\infty$-geometric formula.

Another way of saying this is that a formula is $\infty$-geometric iff it does not contain any of the particles $\to,\neg,\forall,\bigwedge$.
\\[1ex]
 The collection of {\em $\infty$-geometric implications} is generated as follows:
 \begin{enumerate}
 \item If $A,B$ are $\infty$-geometric formulae then $A$, $\neg A$ and $A\to B$ are $\infty$-geometric implications.

  \item If $C(a)$ is a $\infty$-geometric implication and $a$ is a free variable with all occurrences indicated and $x$ does not occur in $C(a)$, then $\forall x\,C(x)$ is
 a $\infty$-geometric implication.

  \item If $\Psi$ is a set of $\infty$-geometric implications having a finite number of free variables then $\bigwedge \Psi$ is a  $\infty$-geometric implication.
\end{enumerate}

 A theory is {\em $\infty$-geometric} if all its axioms are $\infty$-geometric implications.
 }\end{deff}

 \begin{examples}\label{Example-theories}{\em We list some examples of ${\mathcal L}_{\omega_1\omega}$ theories.
 \begin{enumerate}
 \item The theory of {\em torsion groups} is characterized by the  group axioms plus
     the axiom  $$\forall x\,\bigvee\{\underbrace{x\circ\ldots\circ x}_{\mbox{\scriptsize $n$-times}}=e\mid n\geq 1\}.$$
      \item The theory of fields with characteristic a prime is characterized by the  field axioms together with
     the axiom  $$\bigvee\{\underbrace{1+\ldots+1}_{\mbox{\scriptsize $n$-times}}=0\mid n\geq 2\}.$$
     \item The theory of archimedean ordered fields  is characterized by the  ordered field axioms together with
     the axiom  $$\forall x\,\bigvee\{x<\underbrace{1+\ldots +1}_{\mbox{\scriptsize $n$-times}}\mid n\geq 1\}.$$
     \item The class of structures isomorphic to the standard model of Peano arithmetic
     is characterized by the axioms of $\PA$ conjoined with the axiom
     \begin{eqnarray}\label{paunendlich}&&\forall x\,\bigvee\{x=0+\underbrace{1+\ldots +1}_{\mbox{\scriptsize $n$-times}}\mid n\geq 0\}.\end{eqnarray}
     \item The theory of connected graphs has the usual axioms for graphs and additionally has the axiom
         $$\forall x\forall y[x= y \vee \bigvee\{\exists z_0\ldots \exists z_n(x=z_0\wedge y=z_n\wedge Ez_0z_1\wedge\ldots Ez_{n-1}z_n)\mid n\geq 1\}],$$
         where $E$ is a two-place relation such that $Ez_iz_{i+1}$ expresses that there is an edge going from $z_i$ to $z_{i+1}$.

\end{enumerate}
The above theories, with the exception of the fourth example, are $\infty$-geometric. However, in the fourth example the
induction axioms are not really needed as they are implied in infinitary logic by the axiom (\ref{paunendlich}) and the
 axioms of Robinson arithmetic, i.e. the axioms of $\PA$ pertaining to $0,{\mathrm{suc}},+,\cdot$.
 They can be expressed by means of geometric formulae as shown in \ref{3.11}.
}\end{examples}

In the logic ${\mathcal L}_{\infty\omega}$ one has rules for $\bigwedge$ and $\bigvee$
that generalize those for $\wedge$ and $\vee$, respectively. In the sequent calculus version they
can be rendered thus.
\\[3ex]
 $\begin{array}{lcl} \mbox{\red $\bigwedge$-Conjunction} &&
\\ \phantom{A} && \\
\mbox{\AxiomC{$\begin{array}{c} \mbox{\blue $A$},\Gamma\TS
\Delta\;\;\mbox{ and $A\in\Phi$}\end{array}$} \RightLabel{$\bigwedge\,\mbox{L}$}
\UnaryInfC{$\mbox{\red $\bigwedge\Phi$},\Gamma\TS\Delta$}
\DisplayProof}
  & \phantom{A} &
 \mbox{\AxiomC{$\begin{array}{c}
{\Gamma\TS\Delta,\mbox{\blue $A$}\;\;\mbox{ for all $A\in\Phi$}}\end{array}$}
\RightLabel{$\bigwedge\,\mbox{R}$}
\UnaryInfC{$\Gamma\TS\Delta,\mbox{\red $\bigwedge \Phi$}$} \DisplayProof}\\
&\phantom{A}&
\end{array}$
\\[2ex]
 {\red $\bigvee$-Disjunction}
 \\
$\begin{array}{lcl}
\\ \phantom{A} && \\
\mbox{\AxiomC{$\begin{array}{c} \Gamma\TS \Delta,\mbox{\blue
$A$}\;\;\mbox{ and $A\in\Phi$}\end{array}$} \RightLabel{$\bigvee\,\mbox{R}$}
\UnaryInfC{$\Gamma\TS \Delta,\mbox{\red $\bigvee\Phi$}$}
\DisplayProof}
  & \phantom{A} &
 \mbox{\AxiomC{$\begin{array}{c}
{\mbox{\blue $A$},\Gamma\TS\Delta\;\;\mbox{ for all $A\in\Phi$} }\end{array}$}
\RightLabel{$\bigvee\,\mbox{L}$.} \UnaryInfC{$\mbox{\red $\bigvee\Phi$},\Gamma\TS\Delta
$} \DisplayProof}
\end{array}$
\\[2ex]
A detailed proof system for the logic ${\mathcal L}_{\infty\omega}$ will be provided in section 5. Since the technique of cut elimination will be an essential tool
in our investigations, the  sequent calculus is most appropriate.

 \section{Adding the axiom of choice (internally) to geometric theories does not preserve conservativity}
 This section features two examples of geometric theories where the internal addition of $\AC$ does not preserve geometric conservativity.
In the subsequent section we will argue that the external addition of $\AC$, in a certain sense, just amounts to arguing in a stronger background theory.
 It might produce interesting results but perhaps nothing that's not easily obtainable from the Boolean-valued approach to forcing combined with the completeness result for ${\mathcal L}_{\omega_1\omega}$ (both from the 1960s).

\subsection{First example}
 The example to be presented is a first-order theory.
 To define it we draw on a simple method, that is sometimes called {\em Morleyisation},
by which every theory can be given a geometric axiomatization in a richer language.\footnote{One place where one can find this terminology
is Sacks' book from 1972 \cite[p. 256]{sacks}. 
 The technique was used by Skolem in the 1920s and conceivably could have even older roots. Albeit Skolemization would be more appropriate, that name is already used for something else.
Keisler in his 1977 paper \cite[Theorem 2.18]{keisler} refers to this gadget as the introduction of {\em Skolem relations}.  Hodges, in his book \cite[p. 62]{hodges} from 1993,  called this method of gaining a $\forall\exists$ axiomatization and quantifier elimination in a richer language
 {\em atomization}.
For Morleyization in a topos-theoretic setting see e.g. Johnstone's book \cite[p. 858]{johnstone2} from 2002.}

\begin{deff}{\em Let $\mathcal L$ be a language
  In this subsection we shall only be concerned with first order
  formulae. $\forall \vec x\,(A_1(\vec x\,)\rightleftarrows A_2(\vec x\,))$ will stand for two formulae
namely
$\forall \vec x\,(A_1(\vec x\,)\to A_2(\vec x\,))$ and $\forall \vec x\,(A_2(\vec x\,)\to A_1(\vec x\,))$.

 For each formula $A(u_1,\ldots,u_n)$ of $\mathcal L$ with all free variables indicated we add two new $n$-ary relation
 symbols $P_{A(\vec u\,)}$ and $N_{A(\vec u\,)}$ to the language, where $\vec u=u_1,\ldots,u_n$. Call the new language
 $\mathcal L^{a}$. The first-order theory $M^{a}$ in the language $\mathcal L^{a}$ has 
 the following axioms:
 \begin{enumerate} \item $\forall \vec x\,\neg(P_{A(\vec u\,)}(\vec x)\,\wedge\,N_{A(\vec u\,)}(\vec x))$.
 \item $\forall \vec x\,(P_{A(\vec u\,)}(\vec x)\,\vee\,N_{A(\vec u\,)}(\vec x))$.
\item If $A(\vec u\,)$ is atomic add the axioms $\forall \vec x\,(P_{A(\vec u\,)}(\vec x)\rightleftarrows  A(\vec x\,))$.
\item If $A(\vec u\,)$ is $B(\vec u\,)\wedge C(\vec u\,)$ add
 $\forall \vec x\,(P_{A(\vec u\,)}(\vec x\,)\rightleftarrows  P_{B(\vec u\,)}(\vec x\,)\wedge P_{C(\vec u\,)}(\vec x\,))$.

 \item If $A(\vec u\,)$ is $B(\vec u\,)\vee C(\vec u\,)$ add
 $\forall \vec x\,(P_{A(\vec u\,)}(\vec x\,)\rightleftarrows  P_{B(\vec u\,)}(\vec x\,)\vee P_{C(\vec u\,)}(\vec x\,)).$

  \item If $A(\vec u\,)$ is $\neg B(\vec u\,)$ add
 $\forall \vec x\,(P_{A(\vec u\,)}(\vec x\,)\rightleftarrows  N_{B(\vec u)}(\vec x\,)).$

 \item If $A(\vec u\,)$ is $B(\vec u\,)\to C(\vec u\,)$ add
 $\forall \vec x\,(P_{A(\vec u\,)}(\vec x\,)\rightleftarrows   N_{B(\vec u\,)}(\vec x\,)\vee P_{C(\vec u\,)}(\vec x\,)).$

 \item If $A(\vec u\,)$ is $\exists y B(\vec u,y)$ add
 $\forall \vec x\,(P_{A(\vec u\,)}(\vec x\,)\rightleftarrows  \exists y\,P_{B(\vec u,v)}(\vec x,y)).$

 \item If $A(\vec u\,)$ is $\forall  y B(\vec u,y)$ add
 $\forall \vec x\,(N_{A(\vec u\,)}(\vec x\,)\rightleftarrows  \exists y\,N_{B(\vec u,v)}(\vec x,y)).$
\end{enumerate}
If $T$ is a first-order theory, we denote by $T^a$ the theory $M^a$ augmented by the axioms
$$\forall \vec x\,P_{A(\vec u)}(\vec x\,)$$
for all axioms $\forall \vec x\,A(\vec x\,)$ of $T$.

Clearly, $M^a$ and $T^{a}$  are finite geometric theories.
}\end{deff}

\begin{lem}\label{M1} Let $M^a$, $T$ and $T^a$ as above. Let $\vdash^i$ signify intuitionistic deducibility.
\begin{itemize}
\item[(i)] For every formula $A(\vec u\,)$ of $\mathcal L$ with all free variables indicated,
 $$M^{a}\vdash^i \forall \vec x\,[A(\vec x\,) \gdw P_{A(\vec u)}(\vec x\,)].$$
  \item[(ii)] As a classical theory, $T^a$ is conservative over $T$, that is,
 for every $\mathcal L$-sentence $B$,
 $$T\vdash^c B\;\;\mbox{ iff }\;\; T^a\vdash^c B.$$
 This is in general not true for $T$ based on intuitionistic logic.
 \end{itemize}
 \end{lem}
 \prf (i) is proved by induction on the generation of $A(\vec u)$, making use of the excluded middle
 principle for $P_{A(\vec u)}$ that is encapsulated in the first two axioms of $M^a$.

 (ii) This can be shown syntactically but the model-theoretic proof is shorter. Every $\mathcal L$-structure $\mathfrak A$ can be expanded in just one way
 to an $\mathcal L^a$-structure  $\mathfrak A^a$ which is a model of $M^a$.
Hence every model $\mathfrak M$ of $T$ can be expanded in just one way
 to an $\mathcal L^a$-structure. Moreover,  $\mathfrak M^a$  is a model of $T^a$.
 Also, by (i), $T^a$ comprises $T$ as it proves all axioms of $T$.
 \qed

 \begin{cor} Let $\mathcal L$ be the language of set theory and
 $\ZF^a$ be the Morleyization of Zermelo-Fraenkel set theory.
 $\ZF^a+\AC$ is not conservative over $\ZF^a$ for geometric implications
 of $\mathcal L^a$.
 \end{cor}
 \prf Let $\AC$ be the statement $\forall \vec x\,B(\vec x\,)$.
 By Lemma \ref{M1}(i) we have $$\ZF^a\vdash \AC\leftrightarrow \forall \vec x\,P_{B(\vec u\,)}(\vec x\,),$$ and hence
 $\ZF^a+\AC\vdash \forall \vec x\,P_{B(\vec u\,)}(\vec x\,).$
 If $\ZF^a+\AC$ were conservative over $\ZF^a$ for geometric formulae
 we could infer that $\ZF^a\vdash \forall \vec x\,P_{B(\vec u\,)}(\vec x\,)$
 and hence $\ZF^a\vdash\AC$, which would yield
  $\ZF\vdash \AC$ by Lemma \ref{M1}(iv). \qed

  \subsection{Second example} Here we study an infinitary theory.
         Let $\mathcal L'$ be the language with a set of constants $X$ and
         infinitely many unary predicates $P_n$ and $Q_n$ for $n\in{\mathbb N}$.
        Let $T'$ be the $\mathcal L'_{\infty\omega}$-theory
        with the following axioms:
        \begin{itemize}
        \item[(i)] $\forall z\,\neg[P_n(z)\wedge Q_n(z)]$ for all $n\in {\mathbb N}$.
        \item[(ii)] $\bigvee_{a\in X}P_n(a)$  for all $n\in{\mathbb N}$;
        \item[(iii)]  $\bigvee_{n\in \mathbb N} Q_n(f(n))$ for all $f\in X^{\mathbb N}$.
            \end{itemize}
            Note that $f$  does not appear as a function symbol in (iii); $f(n)$ is just a constant from $X$.
 $T'$ is clearly a $\infty$-geometric theory.

In $\mathcal L'_{\infty\omega}$ we can express an instance $A_{cc}$ of countable choice as follows:
        $$\bigwedge_{n\in\mathbb N}\bigvee_{a\in X}P_n(a)\;\to\;\bigvee_{ f\in X^{\mathbb N}}\bigwedge_{n\in\mathbb N} P_n(f(n)).$$
Now observe that $T'+A_{cc}$ is a syntactically inconsistent theory, where the latter means that
an inconsistency $B\wedge\neg B$ for some formula $B$ can be deduced with the help of the usual
logical rules and the infinitary proof rules given at the end of Section \ref{infge}.\footnote{Details of the infinitary proof system will be provided in Section \ref{sequent}.} By contrast,
$T'$ is syntactically consistent as long as $X$ has at least two elements; although $T'$ does not have a model in $\mathbf{Set}$ if we assume countable choice to hold in $\mathbf{Set}$.
That $T$ is syntactically consistent can be seen as follows.
Let $V[G]$ be a forcing extension of the ground model $V$ in which the set $Y:=
X^{\mathbf N}$ of $V$ becomes countable. In $V[G]$ there is an enumeration of all functions $f\in Y$, say $Y=\{f_0,f_1,f_2,\ldots\}$.
Let $g:{\mathbb N}\to X$ be defined in such a way that $g(n)\ne f_n(n)$.
This is possible since $X$ has more than one element.
Now define a model $\mathfrak M$ for $T'$ in $V[G]$ by letting  $M=X$
and interpreting $P_n^{\mathfrak M}$ as $\{g(n)\}$ and $Q_n^{\mathfrak M}$ as $X\setminus\{g(n)\}$. This shows that $T'$ has a model in $V[G]$ and thus $T'$ is syntactically consistent.

\section{Adding the axiom of choice externally to geometric theories does  preserve conservativity}
Let $\mathcal L$ be a language and $T$ be a ${\mathcal L}_{\infty\omega}$-theory. We extend $\mathcal L$ to $\mathcal L'$ by adding
two unary predicate symbols $\mathrm S$ and $\mathrm U$ and the binary relation symbol $\in$. The idea is to define a set theory with urelements where the axioms of $T$ are supposed to hold for the urelements.
 Formally this means that every axiom $A$ of $T$ has to be relativized to $\mathrm U$, denoted $A^{\mathrm U}$,
 i.e. all quantifier occurrences $\forall x \ldots x\ldots $ and $\exists y \ldots y \ldots$ in $A$ have to be replaced by $\forall x({\mathrm U}(x) \to \ldots x\ldots)$ and $\exists x({\mathrm U}(x) \,\wedge\, \ldots x\ldots)$, respectively. $T^{\mathrm U}$ denotes the theory with language $\mathcal L'$
 and all axioms $A^{\mathrm U}$ where $A$ is an axiom of $T$.

  The axioms of set theory then hold for the objects in $S$.
 The axiom of extensionality has to be given in the form
 $$\forall x,y[{\mathrm S}(x)\wedge{\mathrm S}(y)\wedge \forall z(z\in x\leftrightarrow z\in y)\to x=y].$$
 Further axioms proclaim that everything is either an urelement or a set but not both, that urelements have no elements, and that the urelements form a set:
 $\forall x\,[{\mathrm U}(x)\vee {\mathrm S}(x)]$,
 $\forall x\,\neg [{\mathrm U}(x)\wedge {\mathrm S}(x)]$,
 $\forall x,y[{\mathrm U}(x) \to y\notin x]$, $\exists y[{\mathrm S}(y)\wedge\,
 \forall x[x\in y\leftrightarrow {\mathrm U}(x)]$.

 Let $\ZF^{\mathrm U}_n$ denote the set theory with language $\mathcal L'$ with
 urelement axioms having the above axioms, the usual axioms of set theory (Pairing, Union, Foundation, Powerset) expressed for objects of sort $\mathrm S$, Separation extended to
 the language $\mathcal L'$, but with Replacement restricted to $\Sigma_n$-formulae of $\mathcal L'$.

 Below we refer to {\em definable global choice} by which we mean that a formula of  set theory (usually with extra parameters) defines a well-ordering on the entire universe (see \cite[V.3.9]{levy} for details).
 The actual formula will be revealed in the proof of the next theorem.
  We then have the following conservativity result.

 \begin{thm}[$\ZFC$]\label{globalc} Let $B$ be a sentence of $\mathcal L_{\infty\omega}$.
 Then: $$T\vdash B\;\;\; \mbox{ iff } \;\;\;\ZF_n^{\mathrm U}+\mbox{definable global choice}+\GCH+T^{\mathrm U}\vdash B^{\mathrm U}.$$
 If $T$ is $\infty$-geometric and $A$ is a geometric implication, then also
 $$T\vdash^i B \;\;\;\mbox{ iff }\;\;\; \ZF^{\mathrm U}_n+\mbox{definable global choice}+\GCH+T^{\mathrm U}\vdash B^{\mathrm U}.$$
 \end{thm}
 \prf We argue in our background universe satisfying $\ZFC$. Suppose  $\ZF^{\mathrm U}_n+\mbox{global choice}+\GCH+T^{\mathrm U}\vdash B^{\mathrm U}.$ We then switch to a forcing extension $V[G]$ in which the language
 $\mathcal L'$, the formula $B$ and its subformulae as well as the axioms of $T^{\mathrm U}$ together with their subformulae belong a countable transitive set $X$. Let $f:X\to \mathbf N$ be a bijection. Arguing in $V[G]$, we shall work in the relativized constructible hierarchy $L(f)$ which starts with $\mathrm{TC}(f)$, the transitive closure of $\{f\}$ (see \cite[13.24]{jech}). $L(f)$ has a global definable well-ordering since  $\mathrm{TC}(f)$ is countable in $V[G]$. It's also a model of $\GCH$.
 Using the reflection principle of $\ZF$, we can
 take any model $\mathfrak M$ of $T$ in $L(f)$ and expand it into a model of
 $\ZF^{\mathrm U}_n+\mbox{global choice}+\GCH+T^{\mathrm U}$. Thus $\mathfrak M$ will
 satisfy $B$. But in $L(f)$, $T$ is a ${\mathcal L}_{\omega_1\omega}$-theory
 and therefore, by the completeness theorem for this logic,
  there exists a deduction of $B$ from $T$ in $L(f)$. Consequently, if we work in a sequence calculus, invoking Theorem \ref{reductionII} yields that
  there exists also a cut-free deduction of $\bigwedge T\to B$ in $L(f)$,
  where $\bigwedge T$ signifies the conjunction of all axioms of $T$.

  Now it's crucial to observe that $\bigwedge T$ and $B$ both belong to the ground model.
  It remains to show that there is also a deduction of
  $\bigwedge T\to B$ in the ground model $V$. To this end we shall prove a more general result:
  \begin{itemize} \item[$(*)$]
  If there is cut free deduction $\mathcal D$ of the sequent $\Gamma\TS \Delta$ in $V[G]$ and $\Gamma\TS\Delta$ belongs to the ground model $V$, then there already exists a deduction of $\Gamma\TS\Delta$ in $V$.
  \end{itemize}
  In $(*)$ we refer to the sequent calculus for ${\mathcal L}_{\omega_1\omega}$ to be described in
  section 5. A sequent $\Gamma\TS\Delta$ consists of two finite sequences of ${\mathcal L}_{\omega_1\omega}$-formulae $\Gamma$ and $\Delta$.
  We proceed by induction on the rank of $\mathcal D$. It is crucial that $\mathcal D$ contains no cuts
  lest the end sequents of the immediate subderivations of $\mathcal D$ contain formulae that are not in the ground model and the inductive proof breaks down.
  The proof is straightforward except for the cases of a $\bigwedge \mbox{R}$ or $\bigvee \mbox{L}$ inferences
  that require a bit more attention.
  So suppose that the last inference of $\mathcal D$ was $\bigwedge \mbox{R}$. Then $\Delta$ is of the form
  $\Delta_0,\bigwedge\Phi$ and we have deductions $\mathcal D_{A}$ of $\Gamma\TS \Delta_0, A$ for all $A\in \Phi$. With $\Phi\in V$ we also have $A\in V$ for all $A\in \Phi$. Thus inductively for every $A$ in $\Phi$ there exists a deduction $\mathcal D_A'$ of $\Gamma\TS \Delta_0, A$. Using collection and the axiom of choice we can then compose a deduction $\mathcal D'$ of $\Gamma\TS \Delta_0, \bigwedge \Phi$  in $V$.
  The case of an $\bigvee \mbox{L}$ inference is similar.

   The final issue to be resolved is how the forcing extension $V[G]$ can be accessed from the ground model $V$. There are several approaches to this (cf. \cite[Ch. VII.9]{kunen}). One proper formal way is to resort to the Boolean valued approach (cf. \cite{bell}).
Also note that in case the language $\mathcal L$, the theory $T$ and $B$ are all countable, it is not necessary to take a forcing extension. Then $L(f)\subseteq V$ and the main ingredient for proving the theorem is just the completeness of ${\mathcal L}_{\omega_1\omega}$.
 \qed

 \begin{rem}{\em  The declared background theory for the previous Theorem is $\ZFC$, however, $\ZF$ would be sufficient. The axiom of choice can be dropped, though this requires a more careful definition of the notion of infinitary deduction which does not have the axiom of choice built into its very definition. The problem lies with the infinitary rules $\bigwedge\mbox{R}$ and
 $\bigvee \mbox{L}$. $\AC$ is needed when we have to pick exactly one deduction for each of the infinitely many premisses of these inferences. But this can be avoided by allowing non-empty sets of subdeductions of the same end sequent to figure in a deduction. Details will be deferred to section 5.1.
 It should perhaps be mentioned that jettisoning $\AC$ when dealing with infinitary deductions is also important for the Barwise completeness theorem
 (see \cite[III.5]{ba}).}
 \end{rem}

 A possible interpretation of those earlier quotes to the effect that adding the axiom of choice to a geometric theory $T$ does not produce new geometric theorems is that $\AC$ is simply added to an ambient external type theory $\ITT$ which is grafted onto $T$. Here $\ITT$ is the intuitionistic type theory that holds in all toposes
  also known as the {\em internal logic} of toposes (see \cite[II]{lambek-scott}).
  The axiom of choice can be expressed  in the language for
the internal logic (the so-called Mitchell-B\'enabou language) in a straightforward way (\cite[II.6]{lambek-scott}).
If one now assumes that the language of $T$ is incorporated into the Mitchell-B\'enabou language
via relativization to a specific sort ${\mathrm U}$
and one also has an appropriate treatment of the infinite connectives
then one gets the following result.

\begin{cor}
 Let $B$ be a sentence of $\mathcal L_{\infty\omega}$.
 Then:
 \begin{eqnarray*}  \mbox{Internal Logic} + T^{\mathrm U} +\AC \vdash^c B^{\mathrm U}  &\mbox{iff} &\\
  \mbox{Internal Logic} + T^{\mathrm U} \vdash^i B^{\mathrm U} &\mbox{iff } &\\
 T\vdash^i B
\end{eqnarray*}
where $\vdash^c$ and $\vdash^i$ signify classical and intuitionistic derivability, respectively.
\end{cor}
\prf As the internal logic can be interpreted in a small fragment of $\ZF$, this is a consequence of
Theorem \ref{globalc}. \qed

\section{A sequent calculus for ${\mathcal L}_{\infty\omega}$}\label{sequent}
In his thesis Gentzen introduced a form of the sequent calculus and his technique of {\em cut elimination}. The sequent calculus can be generalized to ${\mathcal L}_{\infty\omega}$.

A {\bluff sequent} is an expression {\blue $\Gamma\TS \Delta$} where
 {\blue $\Gamma$} and {\blue $\Delta$} are finite sequences of ${\mathcal L}_{\infty\omega}$-formulae
 {\blue $A_1,\ldots,A_n$} and {\blue $B_1,\ldots, B_m$}, respectively.
We also allow for the possibility that $\Gamma$ or $\Delta$ (or
both) are empty. The empty sequent will be denoted by
$\emptyset$.
 {\blue $\Sigma\TS \Delta$} is read, informally, as $\Gamma$ yields
 $\Delta$ or, rather,
  the {\blue conjunction} of the {\blue $A_i$}
  yields  the {\blue disjunction} of the {\blue $B_j$}.
In particular, we have:
\begin{itemize}
\item If {\blue $\Gamma$} is empty, the sequent asserts the
disjunction of the {\blue $B_j$}. \item If {\blue $\Delta$} is
empty, it asserts the negation of the conjunction of the {\blue
$A_i$}. \item if {\blue $\Gamma$} and {\blue $\Delta$} are both
empty, it asserts the {\blue impossible}, i.e. a {\blue
contradiction}.
\end{itemize}

We use upper case Greek letters {\blue
$\Gamma,\Delta,\Lambda,\Theta,\Xi\ldots$}
  to range over finite sequences of formulae. $\Gamma\subseteq
  \Delta$ means that every formula of $\Gamma$ is also a formula
  of $\Delta$. $\Gamma,A$ stands for the sequence $\Gamma$ extended by $A$.

  Next we list the axioms and rules of the first-order sequent
  calculus.
  \begin{itemize} \item
{\reddish Logical Axioms}
               $$\Gamma,\mbox{\reddish $A$}\TS\Delta, \mbox{\reddish $A$}$$
               where {\blue $A$} is any formula. In point of fact,
               one could limit this axiom to the case of atomic formulae {\blue $A$}.
      \item         {\reddish Cut Rule}
 $$\mbox{ \AxiomC{$\begin{array}{c} {\Gamma\TS
\Delta,\mbox{\blue $A$}}\end{array}$} \AxiomC{$\begin{array}{c}
\mbox{\blue
$A$},\Lambda\TS\Theta\end{array}$}\RightLabel{$\mbox{Cut}$}
  \BinaryInfC{$\Gamma,\Lambda\TS\Delta,\Theta$}
\DisplayProof}$$ 
The formula {\blue $A$} is called the {\blue cut formula} of the
inference. \item {\reddish Structural Rules}
$${\mbox{
\AxiomC{$\Gamma\TS \Delta$} 
  \UnaryInfC{$\Gamma'\TS \Delta'$}
\DisplayProof}} \phantom{AA}\mbox{ if
}\Gamma\subseteq\Gamma',\;\Delta\subseteq\Delta'.$$ A special case
of the structural rule, known as {\em contraction}, occurs when
the lower sequent has fewer occurrences of a formula than the
upper sequent. For instance, $A,\Gamma\TS\Delta,B$ follows
structurally from $A,A,\Gamma\TS\Delta,B,B$.
 \item Rules for Logical Operations \end{itemize}
 $$\begin{array}{ll} 
   \mbox{Left} &\mbox{Right}
\\ \phantom{A} & \\
\mbox{\AxiomC{$\begin{array}{c} {\Gamma\TS \Delta,\mbox{\blue
$A$}}\end{array}$} 
\UnaryInfC{$\mbox{\reddish $\neg A$},\Gamma\TS \Delta$}
\DisplayProof}
  & 
 \mbox{\AxiomC{$\begin{array}{c}
{\mbox{\blue $B$},\Gamma\TS \Delta}\end{array}$}
\UnaryInfC{$\Gamma\TS\Delta,
\mbox{\reddish $\neg  B$}$} \DisplayProof}
\\ \phantom{A} & \\
\mbox{\AxiomC{$\begin{array}{c} {\Gamma\TS \Delta,\mbox{\blue
$A$}}\end{array}$} \AxiomC{$\begin{array}{c} {\mbox{\blue
$B$},\Lambda\TS \Theta}\end{array}$}
\BinaryInfC{$\mbox{\reddish
$A\to B$},\Gamma,\Lambda \TS \Delta,\Theta$} \DisplayProof}
  & 
 \mbox{\AxiomC{$\begin{array}{c}
{\mbox{\blue $A$},\Gamma\TS \Delta,\mbox{\blue $B$}}\end{array}$}
\UnaryInfC{$\Gamma\TS\Delta,
\mbox{\reddish $A\to B$}$} \DisplayProof}
\\ \phantom{A} & \\
\mbox{\AxiomC{$\begin{array}{c} \mbox{\blue $A$},\Gamma\TS
\Delta\end{array}$} 
\UnaryInfC{$\mbox{\reddish $A\wedge B$},\Gamma\TS\Delta$}
\DisplayProof} 
 \mbox{\AxiomC{$\begin{array}{c}
{\mbox{\blue $B$},\Gamma\TS\Delta}\end{array}$}
\UnaryInfC{$\mbox{\reddish $A\wedge B$},\Gamma\TS\Delta$} \DisplayProof}
& 
 \mbox{ \AxiomC{$\begin{array}{c} {\Gamma\TS
\Delta,\mbox{\blue $A$}}\end{array}$} \AxiomC{$\begin{array}{c}
{\Gamma\TS\Delta,\mbox{\blue
$B$}}\end{array}$}
  \BinaryInfC{$\Gamma\TS\Delta,\mbox{\reddish $A\wedge B$}$}
\DisplayProof}
\\ \phantom{A} & \\
\mbox{ \AxiomC{$\begin{array}{c} {\mbox{\blue $A$},\Gamma\TS
\Delta}\end{array}$} \AxiomC{$\begin{array}{c} {\mbox{\blue
$B$},\Gamma\TS\Delta}\end{array}$} 
  \BinaryInfC{$\mbox{\reddish $A\vee B$},\Gamma\TS\Delta$}
\DisplayProof} &

\mbox{\AxiomC{$\begin{array}{c} \Gamma\TS \Delta,\mbox{\blue
$A$}\end{array}$} 
\UnaryInfC{$\Gamma\TS \Delta,\mbox{\reddish $A\vee B$}$}
\DisplayProof}
 \mbox{\AxiomC{$\begin{array}{c}
{\Gamma\TS\Delta, \mbox{\blue $B$}}\end{array}$}
\UnaryInfC{$\Gamma\TS\Delta,
\mbox{\reddish $A\vee B$}$} \DisplayProof}
\\ \phantom{A} & \\

 {\mbox{
\AxiomC{$\begin{array}{c} {\mbox{\blue $F(t)$},\Gamma\TS
\Delta}\end{array}$} \RightLabel{$\forall\,\mbox{\em L}$}
  \UnaryInfC{$ \mbox{\reddish $\forall x\,F(x)$},\Gamma\TS\Delta$}
\DisplayProof}}
 & 
  {\mbox{
\AxiomC{$\begin{array}{c} {\Gamma\TS\Delta, \mbox{\blue
$F(a)$}}\end{array}$} \RightLabel{$\forall\,\mbox{\em R}$}
  \UnaryInfC{$\Gamma\TS\Delta, \mbox{\reddish $\forall x\,F(x)$}$}
\DisplayProof}}
 \\ \phantom{A} & \\
  {\mbox{
\AxiomC{$\begin{array}{c} {\mbox{\blue
$F(a)$},\Gamma\TS\Delta}\end{array}$}
\RightLabel{$\exists\,\mbox{\em L}$}
  \UnaryInfC{$\mbox{\reddish $\exists x\,F(x)$},\Gamma\TS\Delta$}
\DisplayProof}}
 & 
 {\mbox{\AxiomC{$\begin{array}{c}
{\Gamma\TS \Delta,\mbox{\blue $F(t)$}}\end{array}$}
\RightLabel{$\exists\,\mbox{\em R}$}
  \UnaryInfC{$\Gamma\TS\Delta,\mbox{\reddish $\exists x\,F(x)$}$}
\DisplayProof}}

 \end{array}$$
 In {\blue $\forall\mbox{L}$} and {\blue $\exists\mbox{\em R}$}, {\blue $t$} is an
 arbitrary term. The variable {\blue $a$} in
  {\blue $\forall\mbox{R}$} and {\blue $\exists\mbox{\em L}$} is an
  {\reddish eigenvariable} of the respective inference, i.e. {\blue
  $a$} is not to occur in the {\blue lower sequent}.

 The logic ${\mathcal L}_{\infty\omega}$ in addition has rules for $\bigwedge$ and $\bigvee$
that generalize those for $\wedge$ and $\vee$, respectively.
\\[2ex]
 $\begin{array}{lcl} 
\mbox{\AxiomC{$\begin{array}{c} \mbox{\blue $A$},\Gamma\TS
\Delta\;\;\mbox{ and $A\in\Phi$}\end{array}$} \RightLabel{$\bigwedge\,\mbox{L}$}
\UnaryInfC{$\mbox{\red $\bigwedge\Phi$},\Gamma\TS\Delta$}
\DisplayProof}
  & \phantom{A} &
 \mbox{\AxiomC{$\begin{array}{c}
{\Gamma\TS\Delta,\mbox{\blue $A$}\;\;\mbox{ for all $A\in\Phi$}}\end{array}$}
\RightLabel{$\bigwedge\,\mbox{R}$}
\UnaryInfC{$\Gamma\TS\Delta,\mbox{\red $\bigwedge \Phi$}$} \DisplayProof}
\end{array}$
\\[2ex]
$\begin{array}{lcl}
\mbox{\AxiomC{$\begin{array}{c} \Gamma\TS \Delta,\mbox{\blue
$A$}\;\;\mbox{ and $A\in\Phi$}\end{array}$} \RightLabel{$\bigvee\,\mbox{R}$}
\UnaryInfC{$\Gamma\TS \Delta,\mbox{\red $\bigvee\Phi$}$}
\DisplayProof}
  & \phantom{A} &
 \mbox{\AxiomC{$\begin{array}{c}
{\mbox{\blue $A$},\Gamma\TS\Delta\;\;\mbox{ for all $A\in\Phi$} }\end{array}$}
\RightLabel{$\bigvee\,\mbox{L}$} \UnaryInfC{$\mbox{\red $\bigvee\Phi$},\Gamma\TS\Delta
$} \DisplayProof}
\end{array}$
\\[1ex]
In the rules for logical operations,
 the formulae highlighted in the premisses
 are called the {\em minor formulae} of that inference, while the
  formula highlighted in the conclusion is the {\em principal formula} of that
  inference.
 The other formulae of an inference are called {\em side
 formulae}.
 \subsection{What are proofs in ${\mathcal L}_{\infty\omega}$?}
 Proofs in $\mathcal L_{\omega\omega}$ are finite objects and as a result its notion of proof is very robust. For instance, if one knows for a fact that a formula $A$ is provable in a primitive recursive theory $T$,
 then one can conclude that it is inferrable in Heyting arithmetic that $T$ proves $A$.
 In other words, it is immaterial in which background theory (e.g. $\ZFC$ plus large cardinals) we gained the insight that this fact is true. Things are very different when it comes to infinite proofs.
 An example is provided by the ${\mathcal L}_{\infty\omega}$ intuitionistic ($\infty$-geometric) theory ${\mathbf{HA}}_{\infty}$ whose axioms are those of Robinson arithmetic augmented by the axiom $\forall x\,\bigvee_{n\in \mathbb N}x=\bar{n}$, where $\bar{n}$ stands for the $n$-th numeral. In $\ZF$ one can show that there exists an intuitionistic ${\mathbf{HA}}_{\infty}$-proof of a particular statement $A$  that cannot be shown to exist in intuitionistic Zermelo-Fraenkel set theory $\IZF$.\footnote{$A$ can be taken of complexity $\Pi^0_3$, namely $\forall x\exists y\forall z\,(\mathrm{T}(x,x,y)\,\vee\,\neg\mathrm{T}(x,x,z))$, where $\mathrm{T}$ is the predicate from  Kleene's normal form theorem. Now in the presence of countable choice, $A$ implies the existence of a non-computable function. That the existence of a ${\mathbf{HA}}_{\infty}$  proof of $A$
  cannot be shown in $\IZF$ follows from the fact that $\IZF$ plus countable choice is compatible with the statement that all functions from
 $\mathbb N$ to $\mathbb N$ are computable.}

Even in a classical context it may be relevant to choose a suitable formalization of infinite proof.
For instance, for Barwise's completeness theorem for admissible fragments is it important to choose a notion
of proof that ``does not have the axiom of choice built into its very definition" (\cite{ba}, p. 96).

\begin{deff} $(\CZF)$ {\em We will assume that the language $\mathcal L$ is a set. For the formalization of the ${\mathcal L}$-formulae as set-theoretic objects,
we proceed in the same way as Barwise in \cite[III.3]{ba}. They form an inductively defined (proper) class
of sets. An additional assumption we shall make is that the proper subformulae of a formula $A$ are elements of the transitive closure of $A$ which has the pleasant consequence that the rank of a proper subformula of $A$ is
an element of the rank of $A$; in symbols $\rk(B)\in\rk(A)$. Here we use the usual rank definition for sets, i.e.,
$$\rk(a)\,=\,\bigcup \{\rk(x)+1\mid x\in a\}.$$
Note that $\rk(a)$ is always an ordinal and $x+1$ stands for $x\cup\{x\}$. Just as in the classical world, an {\em ordinal} is a transitive set whose elements are transitive. However, the crucial difference between the classical and the intuitionistic context is the forfeiture of the right to use the trichotomy law for ordinals in the latter, i.e., the assertion
$\alpha\in \beta\,\vee\,\beta\in \alpha\,\,\vee\,\alpha=\beta$ can no longer be guaranteed to hold.
}\end{deff}

\begin{deff}{\em  The class of ${\mathcal L}_{\infty\omega}$-proofs (also called {\em deductions} or {\em derivations})  will be defined inductively.
It is desirable to ensure that the inferred formula and the last inference together with its principal and minor formulae are straightforwardly retrievable from any proof $P$.
Firstly, sequents $A_1,\ldots,A_r\Rightarrow B_1,\ldots, B_s$ are easily coded set-theoretically as a pairs of tuples
$\langle \langle A_1,\ldots,A_n\rangle,\langle B_1,\ldots,B_s\rangle\rangle$; we will continue to use
the former notation even when we refer to its set-theoretic coding.

We shall not write down all the clauses for the inductive definition of proofs. Rather we will provide
two illustrative cases, the finitary $(\wedge R)$ and the infinitary $(\bigvee L)$.

Suppose now we have two proofs $\DDD_1$ and $\DDD_2$ of sequents $A_1,\ldots,A_r\Rightarrow B_1,\ldots B_s$
and $A_1,\ldots,A_r\Rightarrow B'_1,\ldots B'_s$ and $1\leq i_0\leq s$  such that $B_i=B'_{i}$ for all $i\ne i_0$. Then
$$\DDD:\;=\;\langle \langle P_1,P_2\rangle, \langle\wedge R,i_0,k_0\rangle, C_1,\ldots,C_p\Rightarrow D_1,\ldots,D_q\rangle$$
is a proof of $ C_1,\ldots,C_p\Rightarrow D_1,\ldots,D_q\rangle$ if
$1\leq k_0\leq q$ and $B_{i_0}\wedge B'_{i_0}=D_{k_0}$, $A_i\in \{C_1,\ldots,C_p\}$ for all $1\leq i\leq r$ and $B_i\in\{D_1,\ldots,D_q\}$ whenever $1\leq i\leq s$ and $i\ne i_0$. It is clear that
from $P$ we can retrieve the last inference together with its principal and minor formulae.

Next assume $\Phi$ is a set of ${\mathcal L}_{\infty\omega}$-formula and there is a function
$f$ with domain $\Phi$ such that there are finite sequences $\Gamma_1,\Gamma_2,\Delta$ of formulae such that for every $A\in \Phi$, $f(A)$ is an inhabited set of proofs
of $\Gamma_1,A,\Gamma_2\Rightarrow \Delta$. Let $i_0$ be the position of any such $A$ in $\Gamma_1,A,\Gamma_2$. Then
$$\DDD'\;:=\;\langle f, \langle \bigvee L,i_0,k_0\rangle, C_1,\ldots,C_p\Rightarrow D_1,\ldots, D_q\rangle$$
is a proof of  $C_1,\ldots,C_p\Rightarrow D_1,\ldots, D_q$ if $1\leq k_0\leq p$, $\bigvee \Phi=C_{k_0}$,
$C\in \{C_1,\ldots,C_p\}$ for every $C\in \Gamma_i$ with $i\in\{0,1\}$, and $D\in \{D_1,\ldots,D_q\}$
for all $D\in \Delta$.

It should by now be obvious how to deal with the other inference  rules.

Observe that the above definition allows to combine each inference step with a structural
rule. This has the advantage that structural rules needn't be treated as separate rules.
There is a lot of leeway as to the details of formalizing infinitary proofs constructively.
However,
observe that the above definition of proof in the case of $(\bigvee L)$ (and dually $(\bigwedge R)$)
contains a crucial part that enables one to construct new proofs from a collection of proofs.
It has the advantage that in $\CZF$ one can prove
from the assumption that there exists a proof of $\Gamma_1,A,\Gamma_2\Rightarrow \Delta$ for every $A\in \Phi$,
that there also exists a proof of $\Gamma_1,\bigvee \Phi,\Gamma_2\Rightarrow \Delta$
by invoking Strong Collection (see \cite{mar,book}). In general, we wouldn't be able to single out a particular proof
of $\Gamma_1,A,\Gamma_2\Rightarrow \Delta$ for every $A\in\Phi$ without relying on the axiom of choice.
\\[1ex]
Let $\vdash \Gamma\Rightarrow \Delta$ signify that there is a proof of $\Gamma\Rightarrow \Delta$.
We also define provability with length $\alpha$ and cut-degree $\rho$, $$\prov{\alpha}{\rho}{\Gamma\Rightarrow \Delta}$$
to mean that there is a proof $P$ of $\Gamma\Rightarrow \Delta$ such that $\rk(P)\in\alpha+1$
and for all cut formulae $C$ in $P$ we have $\rk(C)\in\rho$. In particular  $\prov{\alpha}{0}{\Gamma\Rightarrow \Delta}$ then conveys that there is a proof without cuts.

Naturally, proofs in theories will also be considered. An ${\mathcal L}_{\infty\omega}$-theory $T$ is a set of ${\mathcal L}_{\infty\omega}$-formulae
without free variables. 
The $T$-proofs
are defined as the ${\mathcal L}_{\infty\omega}$-proofs except that one adds the additional axioms  $\Gamma\Rightarrow\Delta, A$ with $A\in T$ to the sequent calculus.

That there is a $T$-proof of $\Phi\Rightarrow \Theta$ will be conveyed by $T\vdash \Phi\Rightarrow\Theta$. Since the theory $T$ can be expressed in ${\mathcal L}_{\infty\omega}$
via a single formula $\bigwedge T$ we also have
\begin{eqnarray}\label{T-als-Satz}
T\vdash\Phi\Rightarrow\Theta &\mbox{ iff }& \vdash \Phi,\bigwedge T\Rightarrow \Theta.
\end{eqnarray}
}\end{deff}

 \section{Turning classical $\infty$-geometric proofs into intuitionistic ones}
 Recall that intuitionistic ${\mathcal L}_{\infty\omega}$-proofs are those obeying the
 simple structural restriction that there can be at most one formula on the right hand side of the sequent symbol $\TS$. Below we shall indicate intuitionistic proofs in ${\mathcal L}_{\infty\omega}$
 by putting $I_{\infty}$ before the turn style symbol.

 The fact that ${\mathcal L}_{\omega\omega}$ geometric proofs can be turned into intuitionistic ones is basically a consequence
 of Gentzen's Hauptsatz. It could have been proved by Gentzen in 1934.
 It is not clear to the present author who first made this observation but it can be found in Orevkov's 1968 paper \cite{orevkov}.\footnote{This is not to say that there are no interesting research questions left.
 Since cut elimination is costly there are still unsolved problems as to how efficient this procedure can be, in general and for special theories (see e.g. \cite{grisha}).}
 As for the ${\mathcal L}_{\infty\omega}$ case it is not clear to him whether there are any syntactic proofs in the published literature (before \cite{prooftheory}). But the purpose of this part of the article, rather than originality, is to show that there is an easy syntactic proof that can also be formalized in the constructive set theory
 $\CZF$ (see \cite{mar,book}). Closer inspection would actually reveal that intuitionistic Kripke-Platek set theory (see \cite{mar,book}) suffices.
 As in the finite case the crucial tool is the Hauptsatz for ${\mathcal L}_{\infty\omega}$.

 \begin{thm}\label{haupt} If $\prov{\alpha}{\rho}{\Gamma\Rightarrow \Delta}$ then there exists $\alpha'$ such that
 $\prov{\alpha'}{0}{\Gamma\Rightarrow \Delta}$.
 \end{thm}
 The proof  of \ref{haupt} in $\CZF$ will be deferred to section \ref{Hauptsatz}.
Without paying attention to constructivity issues, for the countable logic ${\mathcal L}_{\omega_1\omega}$ this was essentially shown by Tait \cite{tait68}.

 The main result of this section requires knowledge of some basic facts.

\begin{lem}[Substitution]\label{substi} Let $\Gamma(a)\TS \Delta(a)$ be a sequent with all occurrences of the free variable $a$ indicated. Let $t$ be an arbitrary term.
If $\provx{}{\alpha}{\rho}{\Gamma(a)\TS \Delta(a)}$
then $\provx{}{\alpha}{\rho}{\Gamma(t)\TS \Delta(t)}$.
\end{lem}
\prf Proceed by induction on $\alpha$.\qed

 \begin{lem}[Inversion]\label{inversion}
  \begin{itemize} \item[(i)] If $\provx{}{\alpha}{\rho}{\Gamma,A\wedge B\TS \Delta}$ then
  $\provx{}{\alpha}{\rho}{\Gamma,A,B\TS \Delta}$.
  \item[(ii)] If $\provx{}{\alpha}{\rho}{\Gamma\TS \Delta,A\wedge B}$ then
  $\provx{}{\alpha}{\rho}{\Gamma\TS \Delta,A}$ and $\provx{}{\alpha}{\rho}{\Gamma\TS \Delta,B}$.
  \item[(iii)] If $\provx{}{\alpha}{\rho}{\Gamma,A\vee B\TS \Delta}$ then
  $\provx{}{\alpha}{\rho}{\Gamma,A\TS \Delta}$ and  $\provx{}{\alpha}{\rho}{\Gamma,B\TS \Delta}$.
  \item[(iv)] If $\provx{}{\alpha}{\rho}{\Gamma\TS \Delta,A\vee B}$ then
  $\provx{}{\alpha}{\rho}{\Gamma\TS \Delta,A,B}$.
   \item[(v)] If $\provx{}{\alpha}{\rho}{\Gamma\TS A\to B,\Delta}$ then
  $\provx{}{\alpha}{\rho}{A,\Gamma\TS \Delta,B}$.
  \item[(vi)] If $\provx{}{\alpha}{\rho}{\Gamma,A\to B\TS \Delta}$ then
  $\provx{}{\alpha}{\rho}{\Gamma\TS \Delta,A}$ and $\provx{}{\alpha}{\rho}{\Gamma,B\TS \Delta}$.
  \item[(vii)] If $\provx{}{\alpha}{\rho}{\Gamma\TS \neg A,\Delta}$ then
  $\provx{}{\alpha}{\rho}{\Gamma,A\TS \Delta}$.
  \item[(viii)] If $\provx{}{\alpha}{\rho}{\Gamma,\neg A\TS \Delta}$ then
  $\provx{}{\alpha}{\rho}{\Gamma\TS \Delta,A}$.
   \item[(ix)] If $\provx{}{\alpha}{\rho}{\Gamma\TS \Delta,\forall x \,B(x)}$ then
  $\provx{}{\alpha}{\rho}{\Gamma\TS \Delta,B(s)}$ for any term $s$.
    \item[(x)] If $\provx{}{\alpha}{\rho}{\Gamma,\exists x\,B(x)\TS \Delta}$ then
    $\provx{}{\alpha}{\rho}{\Gamma,B(s)\TS \Delta}$ for any term $s$.
   \item[(xi)] If $\provx{}{\alpha}{\rho}{\Gamma,\bigvee\Phi\TS \Delta}$ then
    $\provx{}{\alpha}{\rho}{\Gamma,A\TS \Delta}$ for every $A\in \Phi$.
\item[(xii)] If $\provx{}{\alpha}{\rho}{\Gamma\TS \Delta,\bigwedge\Phi}$ then
    $\provx{}{\alpha}{\rho}{\Gamma\TS \Delta,A}$ for every $A\in \Phi$.

  \item[(xiii)] With the exception of (iv), (vi) and (viii) the above inversion properties remain valid
   for the intuitionistic sequent calculus. One half of (vi) also remains valid intutionistically:

 If $\provxii{}{\alpha}{\rho}{\Gamma,A\to B\TSI \Delta}$ then
  $\provxii{}{\alpha}{\rho}{\Gamma,B\TSI \Delta}$.
  \end{itemize}
  \end{lem}
  \prf All can be shown easily by induction on $\alpha$. \qed

Below we use $\bigvee(\Phi,A)$ to stand for $\bigvee (\Phi\cup\{A\})$.
\begin{lem}\label{coherenthilf}
\begin{enumerate}

\item If $\provxii{}{}{}{\Gamma \TS \bigvee(\Phi,F(s))}$, then
$\provxii{}{}{}{\Gamma \TS \bigvee(\Phi,\exists x\,F(x))}$.

    \item If $\provxii{}{}{}{\Gamma\TS \bigvee(\Phi,B)}$ and
    $\provxii{}{}{}{\Gamma\TS \bigvee(\Phi,C)}$, then
     $\provxii{}{}{}{\Gamma\TS \bigvee(\Phi,B\wedge C)}$.

\item If $\provxii{}{}{}{\Gamma\TSI \bigvee(\Phi,A)}$,
then $\provxii{}{}{}{\Gamma, \neg A\TSI \bigvee\Phi}$.

\item If $\provxii{}{}{}{\Gamma,B\TSI \bigvee\Phi}$ and $\provxii{}{}{}{\Gamma\TSI \bigvee(\Phi,A)}$,
then $\provxii{}{}{}{\Gamma,A\to B\TSI \bigvee\Phi}$.

\item If $\provxii{}{}{}{\Gamma\TSI \bigvee(\Phi,A)}$ and $A\in\Theta$,
then $\provxii{}{}{}{\Gamma\TSI \bigvee(\Phi,\bigvee\Theta)}$.

\end{enumerate}
\end{lem}
\prf
(1) We have
\begin{prooftree}
\Axiom$D\fCenter D$
\RightLabel{\scriptsize $(\bigvee R)$}
\UnaryInf$D\fCenter \bigvee(\Phi,\exists x F(x))\;\;\mbox{ for all $D\in\Phi$}$

\Axiom$F(s)\fCenter F(s)$
\RightLabel{\scriptsize $(\exists R)$}
\UnaryInf$F(s)\fCenter\exists x F(x)$
\RightLabel{\scriptsize $(\bigvee R)$}
\UnaryInf$F(s)\fCenter \bigvee(\Phi, \exists x F(x))$

\RightLabel{\scriptsize $(\bigvee L)$}
\BinaryInfC{$\bigvee(\Phi, F(s))\Rightarrow \bigvee(\Phi,\exists x F(x))$}

\end{prooftree} As $\provxii{}{}{}{\Gamma \TS \bigvee(\Phi,F(s))}$, cutting with $\bigvee(\Phi,F(s))$ yields
$\provxii{}{}{}{\Gamma \TS \bigvee(\Phi,\exists x\,F(x))}$.
\\[2ex]
(2) We have
\begin{prooftree}

\Axiom$B, D\fCenter D$
\RightLabel{\scriptsize $(\bigvee R)$}
\UnaryInf$B, D\fCenter \bigvee(\Phi,B\wedge C)\;\;\mbox{ all $D\in\Phi$}$

\Axiom$B,C\fCenter B$
\Axiom$B,C\fCenter C$
\RightLabel{\scriptsize $(\wedge R)$}
\BinaryInf$B,C\fCenter B\wedge C$
\RightLabel{\scriptsize $(\bigvee R)$}
\UnaryInf$B,C\fCenter \bigvee (\Phi,B\wedge C)$

\RightLabel{\scriptsize $(\bigvee L)$}
\BinaryInfC{$B,\bigvee(\Phi, C)\Rightarrow \bigvee(\Phi, B\wedge C)$}

\end{prooftree}
and therefore
\begin{prooftree}

\Axiom$D,\bigvee(\Phi, C)\fCenter D$
\RightLabel{\scriptsize $(\bigvee R)$}
\UnaryInf$D,\bigvee(\Phi, C)\fCenter \bigvee(\Phi,B\wedge C)\;\;\mbox{ all $D\in\Phi$}$

\AxiomC{}
\doubleLine
\UnaryInfC{$B,\bigvee(\Phi, C)\fCenter \bigvee (\Phi,B\wedge C)$}

\RightLabel{\scriptsize $(\bigvee L)$}
\BinaryInfC{$\bigvee(\Phi, B), \bigvee (\Phi,C)\Rightarrow \bigvee(\Phi,B\wedge C)$}

\end{prooftree}
Cuts with $\provxii{}{}{}{\Gamma \TS \bigvee(\Phi, B)}$ and $\provxii{}{}{}{\Gamma \TS \bigvee (\Phi,C)}$ yield the desired outcome
$\provxii{}{}{}{\Gamma\TS \bigvee(\Phi,B\wedge C)}$.
\\[2ex]
(3) is shown as follows: \begin{prooftree}
\Axiom$\Gamma,A\fCenter A$
\RightLabel{\scriptsize $(\neg L)$}
\UnaryInfC{$\Gamma,A,\neg A\fCenter \phantom{A} $}
\RightLabel{\scriptsize $({\mathcal W}_r)$}
\UnaryInfC{$\Gamma,A,\neg A\fCenter \bigvee\Phi$}

\Axiom$\Gamma,B,\neg A\fCenter B$
\RightLabel{\scriptsize $(\bigvee R)$}
\UnaryInfC{$\Gamma,B,\neg A\fCenter \bigvee\Phi\;\;\mbox{ all $B\in\Phi$}$}

\RightLabel{\scriptsize $(\bigvee L)$}
\BinaryInfC{$\Gamma,\bigvee(\Phi,A),\neg A\fCenter \bigvee\Phi$}

\AxiomC{}
\doubleLine
\UnaryInfC{$\Gamma\fCenter \bigvee(\Phi,A)$}

\RightLabel{\scriptsize $(Cut)$}
\BinaryInfC{$\Gamma,\,\neg A\fCenter \bigvee\Phi$}
\end{prooftree}
(4) We have
\begin{prooftree}
\Axiom$\Gamma, A\fCenter A$

\AxiomC{}
\doubleLine
\UnaryInfC{$\Gamma, B\Rightarrow \bigvee\Phi$}

\RightLabel{\scriptsize $(\rightarrow L)$}
\BinaryInfC{$A,\Gamma, A\rightarrow B\Rightarrow \bigvee\Phi$}

\Axiom$\Gamma,C, A\rightarrow B\fCenter C$
\RightLabel{\scriptsize $(\vee R)$}
\UnaryInfC{$\Gamma,C, A\rightarrow B\fCenter \bigvee\Phi\,\;\mbox{ all $C\in\Phi$}$}

\RightLabel{\scriptsize $(\vee L)$}
\BinaryInfC{$\Gamma,A\rightarrow B ,\bigvee(\Phi,A)\Rightarrow \bigvee\Phi$}
\end{prooftree}
Now cutting out $\bigvee(\Phi,A)$ with $\provxii{}{}{}{\Gamma\TSI \bigvee(\Phi,A)}$
yields $\provxii{}{}{}{\Gamma,A\to B\TSI \bigvee\Phi}$.
\\[2ex]
(5) is shown as follows: \begin{prooftree}
\Axiom$\Gamma,A\fCenter A$
\RightLabel{\scriptsize $(\bigvee R)$}
\UnaryInfC{$\Gamma,A\fCenter \bigvee\Theta$}
\RightLabel{\scriptsize $(\bigvee R)$}
\UnaryInfC{$\Gamma,A\fCenter \bigvee(\Phi,\bigvee\Theta)$}

\Axiom$\Gamma,B\fCenter B\;\;\mbox{ ($B\in\Phi$)}$
\RightLabel{\scriptsize $(\bigvee R)$}
\UnaryInfC{$\Gamma,B\fCenter \bigvee(\Phi,\bigvee \Theta)$}

\RightLabel{\scriptsize $(\bigvee L)$}
\BinaryInfC{$\Gamma,\bigvee(\Phi,A)\fCenter \bigvee(\Phi,\bigvee\Theta)$}

\AxiomC{}
\doubleLine
\UnaryInfC{$\Gamma\fCenter \bigvee(\Phi,A)$}

\RightLabel{\scriptsize $(Cut)$}
\BinaryInfC{$\Gamma\fCenter \bigvee(\Phi,\bigvee\Theta)$}
\end{prooftree}
\qed

 \begin{lem}\label{cohe}
 Let  $\Delta$ be a finite set of
 $\infty$-geometric formulae and $\Gamma$ be a finite set of
 $\infty$-geometric implications.
 \\[1ex]
 If $\provx{}{}{}{\Gamma\TS \Delta}$ then $\provxii{}{}{}{\Gamma\TS \bigvee\Delta}$.
 \end{lem}
 \prf Here we rely on Gentzen's Hauptsatz, Theorem \ref{haupt}, for classical ${\mathcal L}_{\infty\omega}$ logic.

 Let $\mathcal D$ be a cut free deduction of $\Gamma\TS \Delta$. The proof proceeds by induction
 on the ordinal height $\alpha$ of ${\mathcal D}$.

 If $\Gamma\TS\Delta$ is an axiom then there exists an atom $A$ such that $A\in \Gamma\cap\Delta$. Thus $\provxii{}{}{}{\Gamma\TS A}$ and therefore, via
 $(\bigvee R)$, we get  to $\provxii{}{}{}{\Gamma\TS \bigvee\Delta}$.

Now suppose that $\Gamma\TS\Delta$ is the result of an inference rule. 
 We inspect the last inference of $\mathcal D$.
 Note that $\forall R,\bigwedge R,\neg R$ and $\to R$ are ruled out since their principal formulae are not $\infty$-geometric
 formulae and would have to occur in the succedent $\Delta$.

 If the last inference was of the form $\forall L$, $\bigwedge L$, $\wedge L$, $\exists L$, $\bigvee L$, or $\vee L$ we can simply apply the
 induction hypothesis to the premisses and re-apply the same inference in the intuitionistic calculus.

 If the last inference was $\exists R$ we apply the induction hypothesis to its premiss  and subsequently use  Lemma \ref{coherenthilf} (1) to
 get the desired result.

 If the last inference was $\wedge R$ we apply
 the induction hypothesis to its premisses and subsequently use  Lemma \ref{coherenthilf} (2).

 If the last inference was $\neg L$ then its minor formula must be $\infty$-geometric. Thus we can apply
 the induction hypothesis to its premiss and subsequently use  Lemma \ref{coherenthilf} (3).

 If the last inference was $\to L$ then apply
 the induction hypothesis to its premisses and subsequently use  Lemma \ref{coherenthilf} (4).

 If the last inference was $\bigvee R$ then apply
 the induction hypothesis to its premisses and subsequently use  Lemma \ref{coherenthilf} (5).

 The case when the last inference was $\vee R$ is similar to the previous one.
 \qed

 \begin{thm}\label{coherentth} Let $T$ be a $\infty$-geometric theory and suppose that there is a classical proof of a
 $\infty$-geometric implication $G$ from $T$.  Then there is an intuitionistic proof of $G$
 from the axioms of $T$.
 \end{thm}
 \prf Below we shall write $T\vdash A$ for $T\vdash \Rightarrow A$ and $T\vdash^i A$ if there is an
 intuitionistic proof of $\Rightarrow A$ from the axioms of $T$.

 We proceed by induction on the buildup of $G$. First suppose that $G$ is of the form $\forall \vec x\,F(\vec x\,)$ where
 $F(\vec a\,)$ is a $\infty$-geometric implication.
 By (\ref{T-als-Satz})
we have $$\provx{}{}{}{\bigwedge T\TS G}$$  where $\bigwedge T$ is
 the conjunction of all axioms of $T$.
 Using the Inversion Lemma \ref{inversion} (ix) we get $\provx{}{}{}{\bigwedge T\TS F(\vec a\,)}$ and hence
 $\provx{T}{}{}{ F(\vec a\,)}.$ The induction hypothesis (since $F(\vec a)$ is a shorter formula
 than $G$) thus yields $T\vdash^i F(\vec a)$ from which $T\vdash^i G$ follows via (several) $\forall R$ inferences.

 Now suppose that $G$ is of the form $\bigwedge\Phi$, where $\Phi$ is a set of $\infty$-geometric
  formulae. By $\bigwedge$-inversion on the right we get
   $$\provx{}{}{}{\bigwedge T\TS H}$$
   for all $H\in \Phi$, and thus inductively we have
   $$\provxii{}{}{}{\bigwedge T\TS H}$$
   for all $H\in\Phi$, so that via $(\bigwedge\,R)$ we arrive at
   $\provxii{}{}{}{\bigwedge T\TS G}$, thus $T\vdash^i G$.

 If $G$ is of the form $\neg G_0$ with $G_0$ $\infty$-geometric we apply the Inversion Lemma \ref{inversion} (vii) to get
 $$\provx{}{}{}{\bigwedge T, G_0\TS}.$$ By Lemma \ref{cohe} we infer that
 $\provxii{}{}{}{\bigwedge T, G_0\TS}$ and thus, by $\neg R$, we have $$\provxii{}{}{}{\bigwedge T\TS \neg  G_0},$$
  thus $T\vdash^i G$.

 If $F$ is of the form $F_0\to F_1$ with $F_i$ geometric formulae
 we apply the Inversion Lemma \ref{inversion} (v) to get
 $$\provx{}{}{}{\bigwedge T, F_0\TS F_1}.$$ By Lemma \ref{cohe} we infer that
 $\provxii{}{}{}{\bigwedge T, F_0\TS F_1}$. By employing $\to R$ we get $\provxii{}{}{}{\bigwedge T\TS  F_0\to F_1}$ and
  hence $T\vdash^i G$.

      \qed

\section{Constructive cut elimination for ${\mathcal L}_{\infty\omega}$}\label{Hauptsatz}
The usual cut elimination proof for ${\mathcal L}_{\infty\omega}$ uses the Veblen functions (see \cite{veb}) $\varphi_{\alpha}$ in order to measure 
 the ``cost" of cut elimination. In a constructive setting, however, one looses the linearity of
ordinals as well as the principle that every inhabited set of ordinals has a least element.
As a result, the definition of analogs of the $\varphi_{\alpha}$ functions has to be carried out in a different way.
A central gadget of cut elimination in infinitary systems is the ``natural" commutative sum of ordinals $\alpha\#
\beta$. Its definition utilizes the Cantor normal form of ordinals to base $\omega$.
This normal form is not available in $\CZF$ (or $\IZF$) and thus a different approach is called for.
We shall have use for the following induction and recursion principle on ordinals, henceforth referred to as {\em $\lhd$-induction} and {\em $\lhd$-recursion}.

\begin{lem} Define $(\alpha,\beta)\lhd(\alpha',\beta')$ by $$\mbox{$\alpha=\alpha'\,\wedge\,\beta\in \beta'$ or
$\alpha\in \alpha'\,\wedge\,\beta=\beta'$ or $\alpha\in\alpha'\,\wedge\,\beta\in\beta'$}.$$
\begin{itemize}
\item[(i)] $(\CZF)$ $\forall \alpha\forall \beta\,
[\forall \gamma\forall \delta((\gamma,\delta)\lhd(\alpha,\beta)\,\to\,F(\gamma,\delta))\;\to\;F(\alpha,\beta)]\,\to\,\forall \alpha\forall \beta\,F(\alpha,\beta)$.
    \item[(ii)] $(\CZF)$ If $G$ is a total $(n+3$-ary class function $G:V^n\times\On\times\On\times V\to\On$
    then there is a (unique) (n+2)-ary class function $F:V^n\times\On\times\On\to \On$
    such that $$F(\vec x,\alpha,\beta)=G(\vec x,\alpha,\beta,\{\langle \gamma,\delta,F(\vec x,\gamma,\delta)\rangle\mid (\gamma,\delta)\lhd(\alpha,\beta)\}).$$
    \end{itemize}
    \end{lem}
\prf (i): Assume
\begin{eqnarray}\label{lhd1}&&\forall \alpha\forall \beta\,
[\forall \gamma\forall \delta((\gamma,\delta)\lhd(\alpha,\beta)\,\to\,F(\gamma,\delta))\;\to\;F(\alpha,\beta)].\end{eqnarray}
Fix an arbitrary ordinal $\rho$. We show \begin{eqnarray}\label{lhd2}&& \forall \xi\in \rho\,F(\alpha,\xi)\end{eqnarray}
by induction on $\alpha\in\rho$. So the inductive assumption gives $\forall \alpha_0\in\alpha
\,\forall \xi\in \rho\,F(\alpha,\xi)$. We then use
use a further subsidiary induction on $\beta\in\rho$ to show $F(\alpha,\beta)$.
 By (\ref{lhd1}) it suffices to show
 \begin{eqnarray}\label{lhd3} && \forall \gamma\forall \delta[(\gamma,\delta)\lhd(\alpha,\beta)\,\to\,F(\gamma,\delta)].\end{eqnarray}
 So suppose that $(\gamma,\delta)\lhd(\alpha,\beta)$.\\
Case 1: $\gamma=\alpha$ and $\delta\in \beta$. $F(\gamma,\delta)$ follows by the subsidiary induction
 hypothesis.\\
Case 2: $\gamma\in\alpha$ and $\delta=\beta$. $F(\gamma,\delta)$ follows by the main induction
 hypothesis.\\
Case 3: $\gamma\in\alpha$ and $\delta\in\beta$. $F(\gamma,\delta)$ also follows by the main induction
 hypothesis.

 Thus we have shown (\ref{lhd3}). This establishes (\ref{lhd2}). Since $\rho$ was arbitrary it follows that
 $F(\alpha,\beta)$ holds for all $\alpha,\beta$.\footnote{The reason for restricting the quantifier in (\ref{lhd2}) to $\rho$ is that it shows that $\lhd$-induction with $F(\alpha,\beta)$ follows from
 $\in$-induction using a formula having no more unbounded quantifiers than $F$. Of course, this is not essential to the current paper.}

 (ii) Noting that $\{(\gamma,\delta)\mid (\gamma,\delta)\lhd(\alpha,\beta)\}$ is a set, (ii) follows from (i) in the same manner as ordinary $\in$-recursion follows from $\in$--induction.
  For more details see \cite{book}. \qed

\begin{deff}{\em  For a class $X$, let $X^{\cup}:=X\,\cup\,\{u\mid \exists y\in X\;u\in y\}$.

Define $\alpha\# \beta$ by $\lhd$-recursion as follows:
 \begin{eqnarray}\label{natural}
 \alpha\#\beta &=& \{\gamma\#\delta\mid (\gamma,\delta)\lhd(\alpha,\beta)\}^{\cup}
 \\ \nonumber &=&
 \{\alpha\# \delta\mid\delta\in\beta\}^{\cup}\,\cup\,\{\eta\#\beta\mid \eta\in\alpha\}^{\cup}\,\cup\,\{\eta\#\delta\mid \eta\in\alpha\,\wedge\,\eta\in\beta\}^{\cup}.\end{eqnarray}}

\end{deff}

\begin{lem} \begin{itemize}
\item[(i)] If $X$ is a set of ordinals then $X^{\cup}$ is an ordinal.
\item[(ii)] $\alpha\# \beta$ is an ordinal and $\alpha\#\beta=\beta\#\alpha$.
\item[(iii)] If $(\gamma,\delta)\lhd (\alpha,\beta)$, then $\gamma\#\delta\in\alpha\#\beta$.
\end{itemize}
\end{lem}
\prf (i) Let $X$ be a set of ordinals. Then $X^{\cup}$ is also a set of ordinals. It remains to show that $X^{\cup}$ is transitive. Suppose $\alpha\in\beta\in X^{\cup}$. Then $\beta\in X$ or $\beta\in \delta$ for some $\delta\in X$. In the first case we have $\alpha\in \bigcup X\subseteq X^{\cup}$. In the second case we
infer that $\alpha\in \delta$ since $\delta$ is an ordinal, thus $\alpha\in \bigcup X\subseteq X^{\cup}$.

(ii) follows by $\lhd$-induction (also using (i)).

(iii) is obvious by definition of $\alpha\#\beta$. \qed

\begin{deff}{\em Let $I$ be a set and $(f_i)_{i\in I}$ be a definable collection of functions $$f_i:\On^{a(i)} \to\On$$
with arity $a(i)\in \mathbb N$. Let $X$ be a set of ordinals. Then the {\em closure of $X$ under $(f_i)_{i\in I}$},
 $\mathrm{Cl}(X,(f_i)_{i\in I})$, is defined as follows:
 \begin{eqnarray*} X_0&= &X^{\cup}\cup \{0\}\\
    X_{n+1} &=& X_n\,\cup\,\{f_i( \alpha_1,\ldots,\alpha_{a(i)})\mid \alpha_1,\ldots \alpha_{a(i)}\in X_n\}^{\cup}\\
    \mathrm{Cl}(X,(f_i)_{i\in I}) &=& \bigcup_{n\in\mathbb N} X_n.
    \end{eqnarray*}
    }\end{deff}

    \begin{lem}\label{cl} Making the same assumptions as in the foregoing definition,
    $\mathrm{Cl}(X,(f_i)_{i\in I})$ is an ordinal which contains $0$ and all elements of $X$.
    Moreover,  $ \mathrm{Cl}(X,(f_i)_{i\in I})$ is closed under $(f_i)_{i\in I}$, i.e., if $\vec \alpha\in  \mathrm{Cl}(X,(f_i)_{i\in I})$
    then   $f_i(\vec \alpha)\in  \mathrm{Cl}(X,(f_i)_{i\in I})$ for all $i\in I$.
    \end{lem}
    \prf Induction on $n$ shows that $X\subseteq X_0\subseteq X_1\subseteq \ldots\subseteq X_n$ and all $X_n$ are ordinals. Hence $\mathrm{Cl}(X,(f_i)_{i\in I})$ is an ordinal.
    If $ \alpha_1,\ldots,\alpha_{a(i)}\in \mathrm{Cl}(X,(f_i)_{i\in I})$, then $\alpha_1,\ldots,\alpha_{a(i)}\in X_n$ for some $n$ since $a(i)\in\mathbb N$, and hence
    $f_i(\alpha_1,\ldots,\alpha_{a(i)})\in X_{n+1}\subseteq \mathrm{Cl}(X,(f_i)_{i\in I})$.
    \qed

    \begin{deff}{\em By main recursion on $\alpha$ and subsidiary recursion on $\beta$ we define the functions
    $$\varphi_{\alpha}:\On\to\On$$ by letting $\varphi_{\alpha}(\beta)$ be the closure of $$\{\varphi_{\alpha}(\xi)\mid \xi\in \beta\}$$
    under the functions $\#$ and $(\varphi_{\eta})_{\eta\in\alpha}$.
    }\end{deff}

    \begin{lem} $(\CZF)$ For all $\alpha,\beta$, $\varphi_{\alpha}(\beta)$ exists.
   Also  $0\in\varphi_{\alpha}(\beta)$.
    \begin{itemize}
    \item[(i)] If $\delta,\xi \in \varphi_{\alpha}(\beta)$ then $\varphi_{\delta}(\xi)\in \varphi_{\alpha}(\beta)$.
        \item[(ii)]  If $\delta\in\beta$   then $\varphi_{\alpha}(\delta)\in \varphi_{\alpha}(\beta)$.

\item[(iii)] If $\delta,\xi\in\varphi_{\alpha}(\beta)$, then $\delta\# \xi\in \varphi_{\alpha}(\beta)$.
    \end{itemize}
    \end{lem}
    The existence of  $\varphi_{\alpha}(\beta)$ follows by main induction on $\alpha$ and subsidiary induction on $\beta$, using Lemma \ref{cl}.
    (i), (ii), and (iii) are immediate by the closure properties of $\varphi_{\alpha}(\beta)$. \qed

\begin{lem}[Reduction] \label{reduction2} $\phantom{A}$ \\ Suppose $\rho=\rk(C)$.
If $\provx{}{\alpha}{\rho}{\Gamma,C\TS \Delta}$ and $\provx{}{\beta}{\rho}{\Xi\TS \Theta,C}$,
then $$\provx{}{\alpha\#\alpha\#\beta\#\beta}{\rho}{\Gamma,\Xi\TS \Delta,\Theta}.$$
\end{lem}
\prf The proof is by induction on $\alpha\#\alpha\#\beta\#\beta$.
We only look at two cases where $C$ was the principal formula of the last
inference in both derivations.

{\bf Case 1}: The first is when $C$ is of the form $\bigwedge\Phi$.
Then we have $$\provx{}{\alpha_1}{\rho}{\Gamma,C,A_0\TS\Delta}$$
and  $$\provx{}{\beta_A}{\rho}{\Xi\TS \Theta,C,A}$$
for some $\alpha_1<\alpha$ and $A_0\in\Phi$  as well as $\beta_A<\beta$
for all  $A\in\Phi$.
By the induction hypothesis we obtain $$\provx{}{\alpha_1\#\alpha_1\#\beta\#\beta}{\rho}{\Gamma,\Xi,A_0\TS \Delta,\Theta}$$
and $$\provx{}{\alpha\#\alpha\#\beta_{A}\#\beta_{A}}{\rho}{\Gamma,\Xi\TS \Delta,\Theta,A_0}.$$
As $A_0$ is a subformula of $C$ we  have $\rk(A_0)\in\rho$. Cutting out $A_0$ thus gives
$\provx{}{\alpha\#\alpha\#\beta\#\beta}{\rho}{\Gamma,\Xi\TS \Delta,\Theta}.$

{\bf Case 2}: The second case is when $C$ is of the form $\forall x\,A(x)$
Then we have $$\provx{}{\alpha_1}{\rho}{\Gamma,C,A(t)\TS\Delta}$$
and  $$\provx{}{\beta_0}{\rho}{\Xi\TS \Theta,C,A(a)}$$
for some $\alpha_1<\alpha$  term $t$ as well as $\beta_0<\beta$
for some eigenvariable $a$. By Lemma \ref{substi}
we have   $$\provx{}{\beta_0}{\rho}{\Xi\TS \Theta,C,A(t)}.$$
By the induction hypothesis we thus get
$$\provx{}{\alpha_1\#\alpha_1\#\beta\#\beta}{\rho}{\Gamma,\Xi,A(t)\TS \Delta,\Theta}$$
and $$\provx{}{\alpha\#\alpha\#\beta_{0}\#\beta_{0}}{\rho}{\Gamma,\Xi\TS \Delta,\Theta,A(t)}.$$
Observing that $\rk(A(t))\in \rk(\forall x\,A(x))=\rho$, cutting out $A(t)$ gives
$$\provx{}{\alpha\#\alpha\#\beta\#\beta}{\rho}{\Gamma,\Xi\TS \Delta,\Theta}.$$
\qed

   \begin{thm}[Cut Elimination Theorem]\label{reductionII} $\phantom{A}$ \\
If $\provx{}{\alpha}{\rho}{\Gamma\TS \Delta}$ then $\provx{}{\vieh_{\rho}(\alpha)}{0}{\Gamma\TS \Delta}$.
   \end{thm}
   \prf We use induction on $\rho$ with a subsidiary induction on $\alpha$.

   If $\Gamma\TS\Delta$ is an axiom then we clearly get the desired result.
   So let's assume that   $\Gamma\TS\Delta$ is not an axiom. Then we have  a last inference
   $(\mathcal I)$ with premisses $\Gamma_i\TS \Delta_i$. Suppose the inference was not a cut.  We then have $\provx{}{\alpha_i}{\rho}{\Gamma_i\TS \Delta_i}$ for some $\alpha_i<\alpha$.
   By the subsidiary induction hypothesis we obtain $\provx{}{\vieh_\rho(\alpha_i)}{0}{\Gamma_i\TS \Delta_i}$.
   Applying the same inference  $(\mathcal I)$ yields $\provx{}{\vieh_{\nu}(\alpha)}{\rho}{\Gamma\TS \Delta}$.

   Now suppose the last inference was a cut with cut formula $C$. Then $\rk(C)\in\rho$ and there exist
   derivations  $\provx{}{\alpha_0}{\rho}{\Gamma_1,C\TS \Delta_1}$
   and  $\provx{}{\alpha_1}{\rho}{\Gamma_2\TS \Delta_2,C}$ for some $\alpha_0,\alpha_1\in\alpha$
   such that $\Gamma_1,\Gamma_2\subseteq \Gamma$ and $\Delta_1,\Delta_2\subseteq \Delta$.
   By the subsidiary induction hypothesis we conclude that
     $\provx{}{\vieh_{\rho}(\alpha_0)}{0}{\Gamma_1,C\TS \Delta_1}$
   and  $\provx{}{\vieh_{\rho}(\alpha_1)}{0}{\Gamma_2\TS \Delta_2,C}$. By the Reduction Lemma \ref{reduction2} we can infer
   $$\provx{}{\vieh_{\rho}(\alpha_0)\#\vieh_{\rho}(\alpha_0)\#\vieh_{\rho}(\alpha_1)
   \#\vieh_{\rho}(\alpha_1)}{\nu}{\Gamma\TS \Delta}$$
   where $\nu=\rk(C)$. Since $\nu\in\rho$ we can now employ the main induction hypothesis,
   yielding  $$\provx{}{\vieh_{\nu}(\vieh_{\rho}(\alpha_0)\#\vieh_{\rho}(\alpha_0)\#\vieh_{\rho}(\alpha_1)
   \#\vieh_{\rho}(\alpha_1))}{0}{\Gamma,C\TS \Delta}.$$
   Since $\vieh_{\nu}(\vieh_{\rho}(\alpha_0)\#\vieh_{\rho}(\alpha_0)\#\vieh_{\rho}(\alpha_1)
   \#\vieh_{\rho}(\alpha_1))\in \vieh_{\rho}(\alpha)$ we arrive at
   $$\provx{}{\vieh_{\rho}(\alpha)}{0}{\Gamma\TS \Delta}.$$
   \qed

\paragraph{Acknowledgement}
Part of the material is based upon research supported by the EPSRC of the UK through grant No. EP/K023128/1. This research was also supported by a Leverhulme Research Fellowship and a Marie Curie International Research Staff Exchange
 Scheme Fellowship within the 7th European Community Framework Programme.
 This publication was made possible through the support of a
grant from the John Templeton Foundation. The opinions expressed in this
publication are those of the author and do not necessarily reflect the
views of the John Templeton Foundation.

Thanks are owed to the anonymous referee of this paper for valuable comments that helped to improve the paper.


\begin{thebibliography}{7}

 \bibitem{mar} P. Aczel, M. Rathjen: {\em Notes on constructive set theory}, Technical Report 40,
Institut Mittag-Leffler (The Royal Swedish Academy of Sciences,Stockholm,2001).
{\em http://www.ml.kva.se/preprints/archive2000-2001.php}

\bibitem{book} P. Aczel, M. Rathjen: {\em Constructive set theory}, book draft,
August 2010.

\bibitem{bah} J. Barwise: {\em An introduction to first-order logic}. In: {\em Handbook of Mathematical Logic} edited by J. Barwise (North-Holland, Amsterdam, 1977) 5--46.

\bibitem{ba} J. Barwise: {\em Admissible Sets and Structures}
(Springer, Berlin 1975).

\bibitem{bell} J.R.  Bell: {\em Boolean-valued models and independence proofs in set theory}.
(Clarendon Press, Oxford, 1977).

\bibitem{carnap} R. Carnap: {\em Formalization of logic}. (Harvard University, 1943)

\bibitem{chang-keisler} C.C. Chang, H.R. Keisler: {\em Model Theory}. $3^{rd}$ edition (North-Holland, Amsterdam, 1990).

\bibitem{gentzen35} G. Gentzen: {\em Untersuchungen \"uber das logische Schliessen I,II.} Mathematische Zeitschrift
39 (1935) 176--210, 405--431.

\bibitem{helmer} O. Helmer: {\em Languages with expressions of infinite length.} Erkenntnis 7 (1938) 138--141.

\bibitem{hodges} W. Hodges: {\em Model theory} (Cambridge University Press, 1993).

\bibitem{jech} T. Jech: {\em Set Theory} $3^{rd}$ edition (Springer, Berlin, 2003)
\bibitem{johnstone} P. Johnstone: {\em Topos theory},
 L.M.S. Monograhs no. 10, (Academic Press, 1977).

\bibitem{johnstone1} P. Johnstone: {\em Sketches of an elephant: A topos theory compendium}, vol. 1 (Clarendon Press, Oxford, 2002).
     p. 899.

 \bibitem{johnstone2} P. Johnstone: {\em Sketches of an elephant: A topos theory compendium}, vol. 2 (Clarendon Press, Oxford, 2002).
     p. 899.



\bibitem{keisler} J. Keisler: {\em Fundamentals of model theory}. In: J. Barwise (ed.): {\em Handbook of Mathematical Logic}
    (North-Holland, Amsterdam, 1977)

    \bibitem{kunen} K. Kunen: {\em Set theory: An introduction to independence proofs}. (North-Holland, Amsterdam, 1980)

    \bibitem{lambek-scott}J. Lambek, P.J. Scott: {\em Introduction to Higher Order Categorical Logic} (Cambridge University Press, 1988).

\bibitem{levy} A. Levy: {\em Basic set theory}. (Springer, Berlin, 1979).


\bibitem{mm} S. Mac Lane, I. Moerdijk: {\em Sheaves in Gemetry and Logic}
(Springer Verlag, 1992).

\bibitem{grisha} G. Mints: {\em Classical and intuitionistic geometric logic}. Talk at Conference on {\em Philosophy, Mathematics, Linguistics: Aspects of Interaction 2012}, 22 May 2012.

\bibitem{moore} G.H. Moore: {\em The prehistory of infinitary logic: 1885 -- 1955}. In: Dalla Chiara, Doets, Mundici, Benthem (eds.)
	{\em Structures and norms in science}. (Kluwer, Dordrecht, 1997).

\bibitem{orevkov} V.P. Orevkov: {\em Glivenko's sequence classes}. In V.P. Orevkov (ed.), { Trudy Mat. Inst. Steklov}, vol. 98, Logical-mathematical caculus. Part 1 (Leningrad, Nauka, 1968) 131--154.

    \bibitem{prooftheory} M. Rathjen: {\em Proof Theory}, Lecture notes for a postgraduate course
    (Leeds, 2012).

\bibitem{sacks} G. Sacks: {\em Saturated model theory}
(Benjamin, 1972)

\bibitem{sch60}K. Sch\"utte: {\em Beweistheorie} (Springer, Berlin,
1960).



\bibitem{tait68} W.W. Tait:
 {\em Normal derivability in classical logic}.
In: J. Barwise: {\em The syntax and semantics of
 infinitary languages}. Lecture Notes in Mathematics 72 (Springer, Berlin,1968) 204--236






\bibitem{veb}O. Veblen: {\em Continuous increasing functions of finite and
transfinite ordinals,} Trans. Amer. Math. Soc. 9 (1908) 280--292.

\bibitem{wraith} G.C. Wraith: {\em Intuitionistic Algebra: Some Recent Developments
in Topos Theory.} Proceedings of the International Congress of Mathematicians
Helsinki, 1978 331--337.




\end{thebibliography}
\end{document}